\documentclass[11pt, article]{amsart}
\usepackage{amsmath,amsfonts,amsthm,amssymb,amscd}
\usepackage[latin1]{inputenc}
\usepackage{mathrsfs}
\usepackage{latexsym}
\usepackage{graphicx}
\usepackage{enumerate}
\usepackage{hyperref}
\usepackage{pgf}
\usepackage{tikz}
\usetikzlibrary{arrows}
\usetikzlibrary{decorations.pathmorphing}
\usetikzlibrary{backgrounds}
\usetikzlibrary{positioning}
\usetikzlibrary{fit}
\usepackage{graphicx}

\def\M{\mathbb{N}_0}
\def\N{\mathbb{N}}

\newcommand{\mex}{\operatorname{mex}} 

\theoremstyle{definition}
\newtheorem{Thm}{Theorem}[section]

\newtheorem{Cor}[Thm]{Corollary}

\newtheorem{Ex}{Example}
\newtheorem{Rem}{Remark}
\newtheorem{Prob}{Problem}

\begin{document}
\title{Comply subtraction games avoiding arithmetic progressions}

\author{Urban Larsson}
\email{urban.larsson@chalmers.se} 
\address{Mathematical Sciences, Chalmers University of Technology 
and University of Gothenburg, G\"oteborg, Sweden}
\date{\today }

\begin{abstract}
Impartial subtraction games on the nonnegative integers have been studied by many and discussed in detail in for example the remarkable work Winning Ways by Conway, Berlekamp and Guy. We describe how \emph{comply} variations of these games, similar to those introduced by Holshouser, Reiter, Smith, St\u anic\u a, can be defined as having its sets of \emph{winning positions} identical to well-known sets \emph{avoiding arithmetic progressions} such as $x+z=2y$, studied by Szerkeres, Erd\H os and Tur\'an, and many others, thus exploring a new territory combining ideas from combinatorial games and combinatorial number theory. The sets we have in mind are \emph{greedy}, that is, for our example: recursively a new nonnegative integer is included to the set if and only if it does not form a three term arithmetic progression with the smaller entries. It is known that the set thus obtained is equivalent to the following log-linear time closed expression: each winning position contains exclusively the digits 0 and 1 in base 3 expansion. In fact this set is impossible as a set of winning positions for a classical subtraction game, in a sense introduced recently by Duch\^ene and Rigo. Therefore our comply-rule generalization of the subtraction games can be seen to resolve new classes of sets as winning positions for heap games. In this context the $\star$-operator for invariant subtraction games was introduced by Larsson, Hegarty and Fraenkel. We define a similar operator for our game. Our comply games generalize into several dimensions. In two dimensions the winning positions can be represented by certain greedy permutations avoiding arithmetic progressions, one of which was recently introduced by Hegarty; while others generalize classical combinatorial games such as Nim and Wythoff Nim.
\end{abstract}

\maketitle
\vskip 30pt

\section{Introduction}\label{s0}
Sets of nonnegative integers avoiding three term arithmetic progressions, that is solutions to $x+z=2y$, have been widely studied by number theorists, but not yet so much by the CGT community. On the other hand, combinatorial number theorists did not yet consider well-known \emph{greedy algorithms} \emph{avoiding} \emph{arithmetic conditions} to any larger extent in the context of \emph{combinatorial games}, such as  for \emph{impartial subtraction games} and the so-called \emph{$\mex$ rule}. In this paper we study interconnections between these fields of mathematics.

We explore generalizations of the classical \emph{subtraction games}---which can be considered as a sub class of the family of so-called \emph{heap games}---under \emph{normal play impartial} rules \cite{BCG82}. We prove that, via their \emph{outcome functions} (a position is in P if and only if none of its options is), they emulate greedily produced sets avoiding arithmetic conditions, such as those avoiding three term arithmetic progressions. The presentation is centered around a version of the so-called ``comply rule" \cite{HoRe01, SmSt02, La11} and (heap-size) invariant games \cite{DuRi10, La12, LHF11} and translation invariant arithmetic conditions; a new $\star$-operator is defined, thus relating our paper again to \cite{LHF11, La12}. For games on two heaps, we demonstrate a very general condition, for a game to have \emph{symmetric} P-positions; that is $(x,y)$ is in P if and only if $(y,x)$ is in P, where $x$ and $y$ denote the respective heap-sizes, represented by nonnegative integers $\M$. These  games are closely related to a certain greedy permutation avoiding ``two dimensional three term arithmetic progressions", originally defined in \cite{He04}. Apart from that, we discuss the general territory, prove some additional results and suggest several problems.

Our games have natural interpretations as ``board games". We begin, by illustrating two examples in Section \ref{s2}. Then in Section \ref{s3} we define our games on one heap of tokens and prove some general results. In Section \ref{s4} we discuss the greedy set avoiding three term arithmetic progressions and how our game emulates it. In Section \ref{s5} we briefly discuss general arithmetic conditions and in Section \ref{s6} we study 2-heap games and greedy injections $\pi:\M\rightarrow \M$ avoiding arithmetic conditions.

\section{Two board games}\label{s2}
Two players, say Black and White, alternate in placing black and white stones in three term arithmetic progressions on equidistance marks representing the nonnegative integers, see Figure \ref{f1}, with ``0" as the left most mark. Three \emph{stones form an arithmetic progression} if and only if the two left most ones have the same color, but different from the single right most stone, and the distance between them equals the distance between the two right most stones. A player in turn moves, by first removing one of the other player's stones and then forming an arithmetic progression with the remaining stone. The top figure shows a legal move option, $8\rightarrow(2,5)$, for White. The second figure from the top illustrates an illegal configuration, namely $8-5\ne 6-2$. However, the first option is not winning for White. In the third figure from the top and downwards we show how white wins, independently of Black's response (the two legal options are illustrated in the 4th and 5th figure from the top). This game is equivalent to the \emph{comply-number game} on three term arithmetic progressions defined in Sections \ref{s3} and \ref{s4}.

\begin{figure} [ht]
\begin{center}
\vspace{0.9 cm}

\begin{tikzpicture}[scale = 0.6]
\filldraw[color = brown, opacity=0.8] (-0.5, -0.1) rectangle (9.5, 1.1);

\draw[black, thin] (0,0.5)--(9,0.5);
\draw[black, thin] (-0.05, 0) -- (-0.05,1); \draw[black, thin] (0, 0) -- (0,1); \draw[black, thin] (1, 0) -- (1,1); \draw[black, thin] (2, 0) -- (2,1); \draw[black, thin] (3, 0) -- (3,1); \draw[black, thin] (4, 0) -- (4,1); \draw[black, thin] (5, 0) -- (5,1); \draw[black, thin] (6, 0) -- (6,1); \draw[black, thin] (7, 0) -- (7,1); \draw[black, thin] (8, 0) -- (8,1); \draw[black, thin] (9, 0) -- (9,1); 
\filldraw[color = yellow, opacity=0.5] (-0.5, -0.1) rectangle (9.5, 1.1);
\foreach \x/\y in {
5/0.5,2/0.5
}
{
\shade[ball color=gray!10!] (\x,\y) circle (13 pt);
\filldraw [color=white, line width=0.3pt, opacity=0.6] (\x,\y) circle (13 pt); 
}

\foreach \x/\y in {
8/0.5
}
\shade[ball color=black!70!] (\x,\y) circle (13 pt);
\end{tikzpicture}

\vspace{0.4 cm}

\begin{tikzpicture}[scale = 0.6]

\draw[black, thin] (0,0.5)--(9,0.5);
\draw[black, thin] (-0.05, 0) -- (-0.05,1); \draw[black, thin] (0, 0) -- (0,1); \draw[black, thin] (1, 0) -- (1,1); \draw[black, thin] (2, 0) -- (2,1); \draw[black, thin] (3, 0) -- (3,1); \draw[black, thin] (4, 0) -- (4,1); \draw[black, thin] (5, 0) -- (5,1); \draw[black, thin] (6, 0) -- (6,1); \draw[black, thin] (7, 0) -- (7,1); \draw[black, thin] (8, 0) -- (8,1); \draw[black, thin] (9, 0) -- (9,1); 
\filldraw[color = gray, opacity=0.75] (-0.5, -0.1) rectangle (9.5, 1.1);
\foreach \x/\y in {
5/0.5,1/0.5
}
{
\shade[ball color=gray!10!] (\x,\y) circle (13 pt);
\filldraw [color=white, line width=0.3pt, opacity=0.6] (\x,\y) circle (13 pt); 
}

\foreach \x/\y in {
8/0.5
}
\shade[ball color=black!70!] (\x,\y) circle (13 pt);
\draw[white] (3,.4) node  {illegal};
\end{tikzpicture}

\vspace{1.3 cm}

\begin{tikzpicture}[scale = 0.6]
\filldraw[color = brown, opacity=0.8] (-0.5, -0.1) rectangle (9.5, 1.1);

\draw[black, thin] (0,0.5)--(9,0.5);
\draw[black, thin] (-0.05, 0) -- (-0.05,1); \draw[black, thin] (0, 0) -- (0,1); \draw[black, thin] (1, 0) -- (1,1); \draw[black, thin] (2, 0) -- (2,1); \draw[black, thin] (3, 0) -- (3,1); \draw[black, thin] (4, 0) -- (4,1); \draw[black, thin] (5, 0) -- (5,1); \draw[black, thin] (6, 0) -- (6,1); \draw[black, thin] (7, 0) -- (7,1); \draw[black, thin] (8, 0) -- (8,1); \draw[black, thin] (9, 0) -- (9,1); 
\filldraw[color = yellow, opacity=0.5] (-0.5, -0.1) rectangle (9.5, 1.1);
\foreach \x/\y in {
4/0.5,0/0.5
}
{
\shade[ball color=gray!10!] (\x,\y) circle (13 pt);
\filldraw [color=white, line width=0.3pt, opacity=0.6] (\x,\y) circle (13 pt); 
}

\foreach \x/\y in {
8/0.5
}
\shade[ball color=black!70!] (\x,\y) circle (13 pt);
\end{tikzpicture}

\vspace{0.9 cm}

\begin{tikzpicture}[scale = 0.6]
\filldraw[color = brown, opacity=0.8] (-0.5, -0.1) rectangle (9.5, 1.1);

\draw[black, thin] (0,0.5)--(9,0.5);
\draw[black, thin] (-0.05, 0) -- (-0.05,1); \draw[black, thin] (0, 0) -- (0,1); \draw[black, thin] (1, 0) -- (1,1); \draw[black, thin] (2, 0) -- (2,1); \draw[black, thin] (3, 0) -- (3,1); \draw[black, thin] (4, 0) -- (4,1); \draw[black, thin] (5, 0) -- (5,1); \draw[black, thin] (6, 0) -- (6,1); \draw[black, thin] (7, 0) -- (7,1); \draw[black, thin] (8, 0) -- (8,1); \draw[black, thin] (9, 0) -- (9,1); 
\filldraw[color = yellow, opacity=0.5] (-0.5, -0.1) rectangle (9.5, 1.1);
\foreach \x/\y in {
4/0.5
}
{
\shade[ball color=gray!10!] (\x,\y) circle (13 pt);
\filldraw [color=white, line width=0.3pt, opacity=0.6] (\x,\y) circle (13 pt); 
}

\foreach \x/\y in {
3/0.5, 2/0.5
}
\shade[ball color=black!70!] (\x,\y) circle (13 pt);
\end{tikzpicture}

\vspace{0.4 cm}

\begin{tikzpicture}[scale = 0.6]
\filldraw[color = brown, opacity=0.8] (-0.5, -0.1) rectangle (9.5, 1.1);

\draw[black, thin] (0,0.5)--(9,0.5);
\draw[black, thin] (-0.05, 0) -- (-0.05,1); \draw[black, thin] (0, 0) -- (0,1); \draw[black, thin] (1, 0) -- (1,1); \draw[black, thin] (2, 0) -- (2,1); \draw[black, thin] (3, 0) -- (3,1); \draw[black, thin] (4, 0) -- (4,1); \draw[black, thin] (5, 0) -- (5,1); \draw[black, thin] (6, 0) -- (6,1); \draw[black, thin] (7, 0) -- (7,1); \draw[black, thin] (8, 0) -- (8,1); \draw[black, thin] (9, 0) -- (9,1); 
\filldraw[color = yellow, opacity=0.5] (-0.5, -0.1) rectangle (9.5, 1.1);
\foreach \x/\y in {
4/0.5
}
{
\shade[ball color=gray!10!] (\x,\y) circle (13 pt);
\filldraw [color=white, line width=0.3pt, opacity=0.6] (\x,\y) circle (13 pt); 
}

\foreach \x/\y in {
0/0.5, 2/0.5
}
\shade[ball color=black!70!] (\x,\y) circle (13 pt);
\end{tikzpicture}

\vspace{0.9 cm}

\begin{tikzpicture}[scale = 0.6]
\filldraw[color = brown, opacity=0.8] (-0.5, -0.1) rectangle (9.5, 1.1);

\draw[black, thin] (0,0.5)--(9,0.5);
\draw[black, thin] (-0.05, 0) -- (-0.05,1); \draw[black, thin] (0, 0) -- (0,1); \draw[black, thin] (1, 0) -- (1,1); \draw[black, thin] (2, 0) -- (2,1); \draw[black, thin] (3, 0) -- (3,1); \draw[black, thin] (4, 0) -- (4,1); \draw[black, thin] (5, 0) -- (5,1); \draw[black, thin] (6, 0) -- (6,1); \draw[black, thin] (7, 0) -- (7,1); \draw[black, thin] (8, 0) -- (8,1); \draw[black, thin] (9, 0) -- (9,1); 
\filldraw[color = yellow, opacity=0.5] (-0.5, -0.1) rectangle (9.5, 1.1);
\foreach \x/\y in {
0/0.5, 1/0.5
}
{
\shade[ball color=gray!10!] (\x,\y) circle (13 pt);
\filldraw [color=white, line width=0.3pt, opacity=0.6] (\x,\y) circle (13 pt); 
}

\foreach \x/\y in {
2/0.5
}
\shade[ball color=black!70!] (\x,\y) circle (13 pt);
\end{tikzpicture}

\end{center}\caption{A board game on three term arithmetic progressions}
\end{figure}
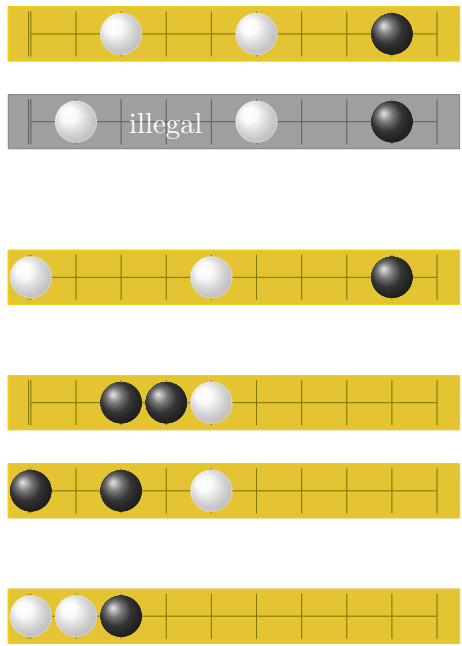\label{f1}

Our second example is the game of \emph{line-Nim}, illustrated in Figures \ref{f2} and~\ref{f3}. It is also discussed in Section \ref{s6} (the max-variation) and in particular in Example~\ref{ex4}. The rules are similar to those of the three term arithmetic progression game, but here the players are allowed to place two stones on a lattice point (a crossing of a horizontal and a vertical line) strictly to the left and below of the remaining stone if and only if all three stones lie on a straight line, all stones must remain in the first quadrant at all times of course. Nim type moves are also allowed, that is, decrease precisely one of the coordinates (see also Example~\ref{e1}). 

The two game boards to the left in Figure \ref{f2} represent move options for White from position $(6,8)$ in the game of \emph{line-Nim}, namely to $\{(3,2),(5,6)\}$ and $\{(6,3)\}$, representing a line-type and Nim type move, respectively. The lower left corner represents position $(0,0)$. To the right we exemplify an illegal configuration. Namely, the black stone is not centered on the line defined by (the centers of) the white stones. Say that White choses the left most move from Figure \ref{f2}. Then we are in Figure \ref{f3}. Black removes one of the white stones and defines the line $(2,3),(4,5)$ which intersects the remaining black stone at $(5,6)$. From this position there are in total 23 options, but in fact (as we will also see in Example \ref{e1} on page \pageref{e1}) the move is winning for Black. 
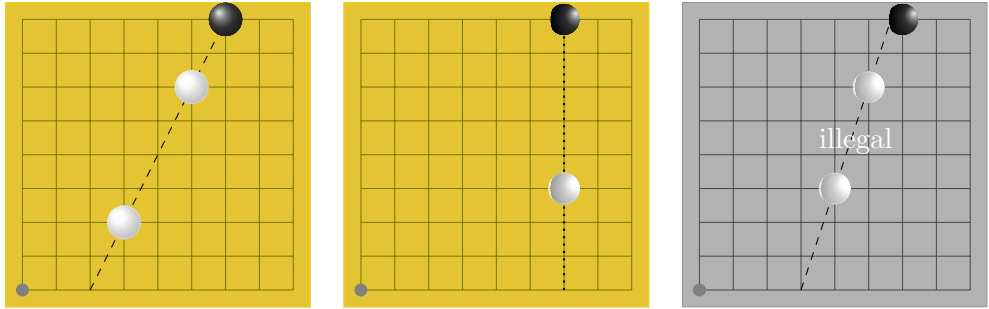
\begin{figure} [h]
\begin{center}

\vspace{0.9 cm}

\begin{tikzpicture}[scale = 0.45]
\begin{scope}[xshift=0cm, yshift=-0cm]
\filldraw[color = brown, opacity=0.8] (-0.5, -0.5) rectangle (8.5, 8.5);

\draw[step=1cm, black] (0, 0) grid (8,8); 
\filldraw[color = yellow, opacity=0.5] (-0.5, -0.5) rectangle (8.5, 8.5);
\filldraw[color=gray] (0, 0) circle (5 pt);

\draw[dashed] (2,0)--(6,8);

\foreach \x/\y in {
3/2, 5/6
}
{
\shade[ball color=gray!10!] (\x,\y) circle (14 pt);
\filldraw [color=white, line width=0.3pt, opacity=0.6] (\x,\y) circle (14 pt); 
}

\foreach \x/\y in {
6/8
}
\shade[ball color=black!70!] (\x,\y) circle (14 pt);
\hspace{1 cm}
\end{scope}
\begin{scope}[xshift=10cm, yshift=-0cm]
\filldraw[color = brown, opacity=0.8](-0.5, -0.5) rectangle (8.5, 8.5);

\draw[step=1cm, black] (0, 0) grid (8,8); 
\filldraw[color = yellow, opacity=0.5] (-0.5, -0.5) rectangle (8.5, 8.5);
\filldraw[color=gray] (0, 0) circle (5 pt);

\draw[dotted, color=black, thick] (6,0)--(6,8);
\foreach \x/\y in {
6/3
}
{
\shade[ball color=gray!10!] (\x,\y) circle (13 pt);
\filldraw [color=white, line width=0.3pt, opacity=0.6] (\x,\y) circle (13 pt); 
}

\foreach \x/\y in {
6/8
}
\shade[ball color=black!70!] (\x,\y) circle (13 pt);
\end{scope}
\begin{scope}[xshift=20cm, yshift=-0cm]
\draw[step=1cm, black, thin] (0, 0) grid (8,8); 
\filldraw[color = gray, opacity=0.6] (-0.5, -0.5) rectangle (8.5, 8.5);
\filldraw[color=gray] (0, 0) circle (5 pt);

\draw[dashed] (3,0)--(5.66,8);
\foreach \x/\y in {
4/3, 5/6
}
{
\shade[ball color=gray!0!] (\x,\y) circle (13 pt);
\filldraw [color=white, line width=0.3pt, opacity=0.65] (\x,\y) circle (13 pt); 
}

\foreach \x/\y in {
6/8
}
\shade[ball color=black!70!] (\x,\y) circle (13 pt);
\draw[white] (2.4,4.4) node  {illegal};
\end{scope}

\end{tikzpicture}

\end{center}\caption{Two move options of the line-Nim game and (to the right) an illegal configuration.}

\end{figure}\label{f2}

\begin{figure} [h]
\begin{center}

\vspace{1.4 cm}

\begin{tikzpicture}[scale = 0.45]
\begin{scope}[xshift=0cm, yshift=-0cm]
\filldraw[color = brown, opacity=0.8](-0.5, -0.5) rectangle (8.5, 8.5);

\draw[step=1cm, black] (0, 0) grid (8,8); 
\filldraw[color = yellow, opacity=0.5] (-0.5, -0.5) rectangle (8.5, 8.5);
\filldraw[color=gray] (0, 0) circle (5 pt);

\foreach \x/\y in {
3/2, 5/6
}
{
\shade[ball color=gray!10!] (\x,\y) circle (14 pt);
\filldraw [color=white, line width=0.3pt, opacity=0.6] (\x,\y) circle (14 pt); 
}
\end{scope}

\begin{scope}[xshift=10cm, yshift=-0cm]
\filldraw[color = brown, opacity=0.8] (-0.5, -0.5) rectangle (8.5, 8.5);

\draw[step=1cm, black] (0, 0) grid (8,8); 
\filldraw[color = yellow, opacity=0.5] (-0.5, -0.5) rectangle (8.5, 8.5);
\filldraw[color=gray] (0, 0) circle (5 pt);

\foreach \x/\y in {
5/6
}
{
\shade[ball color=gray!10!] (\x,\y) circle (13 pt);
\filldraw [color=white, line width=0.3pt, opacity=0.6] (\x,\y) circle (13 pt); 
}

\foreach \x/\y in {
4/5, 2/3
}
\shade[ball color=black!70!] (\x,\y) circle (13 pt);

\end{scope}

\begin{scope}[xshift=20cm, yshift=-0cm]
\filldraw[color = brown, opacity=0.8] (-0.5, -0.5) rectangle (8.5, 8.5);

\draw[step=1cm, black] (0, 0) grid (8,8); 
\filldraw[color = yellow, opacity=0.5] (-0.5, -0.5) rectangle (8.5, 8.5);
\filldraw[color=gray] (0, 0) circle (5 pt);

\foreach \x/\y in {
2/3,1/2
}
{
\shade[ball color=gray!10!] (\x,\y) circle (13 pt);
\filldraw [color=white, line width=0.3pt, opacity=0.6] (\x,\y) circle (13 pt); 
}

\foreach \x/\y in {
4/5
}
\shade[ball color=black!70!] (\x,\y) circle (13 pt);

\end{scope}

\end{tikzpicture}

\end{center}\caption{A game of line-Nim}\label{f3}

\end{figure}
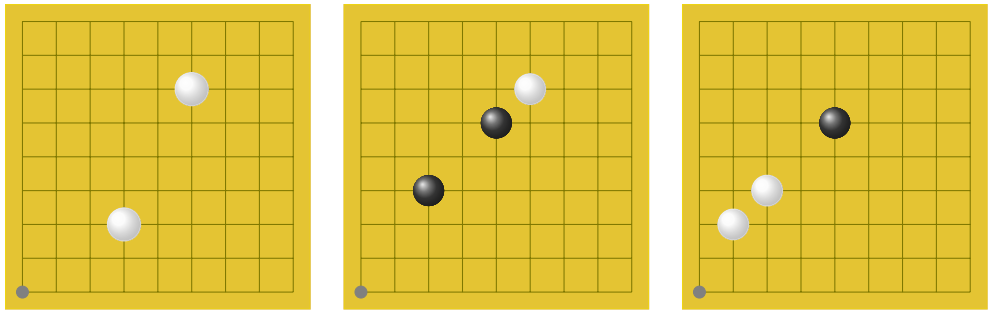

\begin{figure} [h]
\begin{center}

\vspace{0.4 cm}

\begin{tikzpicture}[scale = 0.45]
  
\begin{scope}[xshift=0cm, yshift=-0cm]
\filldraw[color = brown, opacity=0.8](-0.5, -0.5) rectangle (8.5, 8.5);

\draw[step=1cm, black] (0, 0) grid (8,8); 
\filldraw[color = yellow, opacity=0.5] (-0.5, -0.5) rectangle (8.5, 8.5);
\filldraw[color=gray] (0, 0) circle (5 pt);

\foreach \x/\y in {
1/2
}
{
\shade[ball color=gray!10!] (\x,\y) circle (13 pt);
\filldraw [color=white, line width=0.3pt, opacity=0.6] (\x,\y) circle (13 pt); 
}

\foreach \x/\y in {
1/1
}
\shade[ball color=black!70!] (\x,\y) circle (13 pt);

\end{scope}

\begin{scope}[xshift=10cm, yshift=-0cm]
\filldraw[color = brown, opacity=0.8] (-0.5, -0.5) rectangle (8.5, 8.5);

\draw[step=1cm, black] (0, 0) grid (8,8); 
\filldraw[color = yellow, opacity=0.5] (-0.5, -0.5) rectangle (8.5, 8.5);
\filldraw[color=gray] (0, 0) circle (5 pt);

\foreach \x/\y in {
0/1
}
{
\shade[ball color=gray!10!] (\x,\y) circle (13 pt);
\filldraw [color=white, line width=0.3pt, opacity=0.6] (\x,\y) circle (13 pt); 
}

\foreach \x/\y in {
1/1
}
\shade[ball color=black!70!] (\x,\y) circle (13 pt);

\end{scope}
\begin{scope}[xshift=20cm, yshift=-0cm]
\filldraw[color = brown, opacity=0.8] (-0.5, -0.5) rectangle (8.5, 8.5);

\draw[step=1cm, black] (0, 0) grid (8,8); 
\filldraw[color = yellow, opacity=0.5] (-0.5, -0.5) rectangle (8.5, 8.5);
\filldraw[color=gray] (0, 0) circle (5 pt);

\foreach \x/\y in {
0/1
}
{
\shade[ball color=gray!10!] (\x,\y) circle (13 pt);
\filldraw [color=white, line width=0.3pt, opacity=0.6] (\x,\y) circle (13 pt); 
}

\foreach \x/\y in {
0/0
}
\shade[ball color=black!70!] (\x,\y) circle (13 pt);

\end{scope}

\end{tikzpicture}

\end{center}\caption{in which Black wins}\label{f3}

\end{figure}
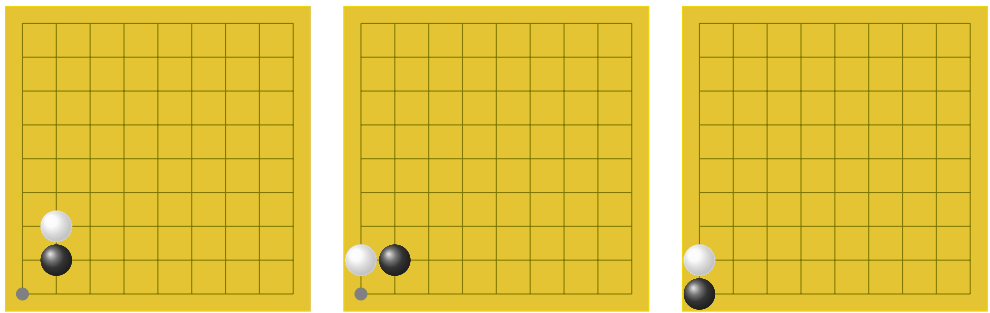

\clearpage

\section{Comply subtraction games}\label{s3}
The rules of a \emph{subtraction game}, here a \emph{game on one heap of tokens}, are as follows: let $S$ denote a set of positive integers and $x$ a nonnegative integer. Then a move option from a heap of $x$ tokens is to a heap of $x-s\ge 0$ tokens, for some $s\in S$. A heap-size is in P (a previous player win) if the player to move cannot win given best play, otherwise it is in N (a next player win). Thus the empty heap is in P and in general the \emph{nim-values} for the positions will be obtained by a \emph{minimal exclusive} algorithm in the following way. Heap-size $x$ has nim-value $g(x)=\text{mex}\{g(x-s)\mid x-s \ge 0, s\in S\}$, where $\text{mex}X=\min(\M\setminus X)$. It follows that heap-size $x$ is in P if and only if $g(x)=0$. Since a position is in P if and only if each option is in N, the \emph{outcome class} (N or P) for an impartial game can always be determined recursively without computing the general nim-values of the N-positions. The following observations are worth mentioning.
\begin{Thm}\label{t1}
Let $A\subset \M$. The following statements are equivalent:
\begin{itemize}
\item [(i)] There is an $a\not\in A$ such that for all $b\in A$ with $b < a$ there is a pair of integers $x,y \in A$ such that $y-x=a-b$. 
\item [(ii)] There is no subtraction game with $A$ as its set of P-positions.
\end{itemize}
\end{Thm}
\noindent{\bf Proof.} If (i) holds then we cannot define a move from $a$ to any position in $A$ since then there would also be a move from $x$ to $y$ which implies that $x$ or $y$ is in N. 

For the other direction, suppose that for all $a\not\in A$ there is a $b\in A$ such that there is no pair of integers $x,y \in A$ with $y-x=a-b>0$. Then, for all $a\not\in A$ we can define $S$ by letting $a-b\in S$.
\hfill $\Box$\\

A corresponding classification for candidate nim-values follows.

\begin{Thm}\label{t2}
Let $f:\M\rightarrow \M$. The following statements are equivalent:
\begin{itemize}
\item [(i)] There is an $a$ with $f(a) > 0$ such that for each set $\{b<a\}$ satisfying mex$\{f(b)\}=f(a)$ there is some $f(b)$ such that for all $b'$ satisfying $f(b') = f(b)$ there is a pair of heap-sizes $x, y$ with $f(x) = f(y) = 0$ such that $x-y=a-b'$. 
\item [(ii)] There is no subtraction game with $f$ as its nim-value function.
\end{itemize}
\end{Thm}
\noindent{\bf Proof.} If (i) holds and f were the nim-value function and there was a move from $a$ to any position with nim-value $g(b)=f(b)<g(a)$ then there would also be a move from $x$ to $y$ which implies that $x$ or $y$ is in N, which contradicts their definition. 

For the other direction, suppose that for all $a$ with $f(a) > 0$ there is a set $\{b<a\}$ satisfying mex$\{f(b)\}=f(a)$ and for all $f(b)<f(a)$ there is at least one $b'$ with $f(b') = f(b)$ such that for all $x, y$ with $f(x) = f(y) = 0$, we have that $x-y\ne a-b'$. Then for all such pairs $a,b'$ we can let $a-b'\in S$. \hfill $\Box$\\

Let us generalize the subtraction games. Let $S$ denote a family of finite subsets of $\N$ and let the heap-size be $x\in \M$. Then a move consists in two parts. The player to move (next player) proposes a set $s\in S$ satisfying 
\begin{align}\label{xs}
x\ge \max s. 
\end{align}
and the player not to move (previous player) chooses one of the numbers $s_i\in s$ to subtract from the given heap-size $x$. Thus the next heap has $x-s_i\ge 0$ tokens. The game ends when the next player cannot propose a set $s$ satisfying (\ref{xs}). The game can equivalently be interpreted as a blocking variation e.g. \cite{HoRe01, SmSt02}, where at each stage of game the player in turn suggests a set $s\in S$ and the other player blocks off all except one $s_i$ from $s$. Note that the set $S$ defines the game and remains fixed, but each proposal of move set $s\in S$ is forgotten when both parts of a move is carried out. We call this the \emph{comply-number game}. 
\begin{Thm}\label{t3}
Suppose that the heap in the comply-number game defined by $S\subset 2^\N$ has $x\in \M$ tokens. Then $x$ is in N if and only if there is some $s\in S$ such that all $x-s_i$ are in P.
\end{Thm}
\noindent{\bf Proof.} The player in turn chooses a set $s\in S$ if applicable and the other player tries to find an N-position of the form $x-s_i$ for some $s_i\in s$. 
\hfill $\Box$\\
\begin{Rem}
A subtraction game is a comply-number game where each move set contains precisely one positive integer.
\end{Rem}
The \emph{dual} of the comply-number game is the \emph{comply-set game}. It is defined accordingly. Let $S$ denote a family of finite subsets of $\N$ and let the heap-size be $x\in \M$. Then a move consists in two parts. The previous player proposes a set $s\in S$ satisfying (\ref{xs}) and the next player chooses one of the numbers $s_i\in s$ to subtract from $x$. Thus the next heap has $x-s_i\ge 0$ tokens. The game ends when the previous player cannot propose a set $s$ satisfying (\ref{xs}). In analogy to the comply-number game this part can equivalently be interpreted as a blocking variation, but the comply interpretation is more direct and gives somewhat simpler game rules. Note that the set $S$ defines the game and remains fixed, but each proposal of move set $s\in S$ is forgotten when both parts of a move is carried out. 
\begin{Thm}\label{t4}
Suppose that the heap in the comply-set game defined by $S\subset 2^\N$ has $x\in \M$ tokens. Then $x$ is in N if and only if for all $s\in S$ with $x\ge\max s$ there is an $s_i$ such that $x-s_i$ is in P.
\end{Thm}
\noindent{\bf Proof.} The player not in turn chooses a set $s\in S$ if applicable and the player in turn tries to find a P-position of the form $x-s_i$ for some $s_i\in s$. \hfill $\Box$\\

By Theorem \ref{t3}, an algorithm to compute the outcomes for the \emph{comply-number game} is as follows: suppose that the outcomes are known for all heap-sizes $x<n$. Then $n$ is in P if and only if for all sets $s\in S$ with $n\ge\max s$ there is a heap-size $n-s_i$ in N. If we wish to compute the outcomes for the \emph{comply-set game} we can use the same algorithm, except exchange P for N and N for P, that is, by Theorem \ref{t4}: suppose that the outcomes are known for all heap-sizes $x<n$. Then $n$ is in N if and only if for all sets $s\in S$ with $n\ge\max s$ there is a heap-size $n-s_i$ in P. We get the following result.

\begin{Thm}\label{PN}
The heap-size $x$ is in P for the comply-set game if and only it is in N for the comply-number game.
\end{Thm}

For a similar result of correspondence between N and P positions for comply- versus blocking impartial games, see \cite{La11}.
\section{Games and sets avoiding arithmetic conditions}\label{s4}
Consider the following general problem \cite{DuRi10,LHF11,La12}.

\begin{Prob}
Let $A\subset \M$. Is there a non-trivial normal play impartial heap-game with $A$ as its set of P-positions?
\end{Prob}

A trivial heap-game is easy to find if $0\in A$, by for example letting each heap-size $\not\in A$ move to $0$. However, a nice property for a heap-game is that each move option is independent of from which heap-size moved from (provided that the move results in a nonnegative heap-size). The main interest for our games is this property of (\emph{heap-size}) \emph{invariance}. In particular it is satisfied by the games we have discussed in Section \ref{s3}. 

Let $A = \{0,1,3,4,9,10,12,13,\ldots\}$. This set corresponds to a known greedy construction which produces a set $A$ of nonnegative integers avoiding \emph{ three-term arithmetic progressions}, that is not containing triples of non-negative integers of the form $(x,x+d,x+2d)$ with $d>0$ the \emph{discrepancy} of the progression. Then $\#(A\cap \{0,\ldots , n\}) \ll n^{\log 2/\log 3}$. Namely this set consists of all integers with digits 0 or 1 in base 3 expansion so that local maxima are obtained at $\frac{3^t-1}{2}$, for all nonnegative integers $t$, $\#A\cap\{0,\ldots ,\frac{3^t-1}{2}\} = 2^t$. 

We look for a normal play impartial heap game that satisfies $$x\in \N\setminus A = \{2,5,6,7,8,11,\ldots\}$$ implies that there is an $y\in A$ such that $x\rightarrow y$ is a legal move, but for all $x,y\in A$ there is no move from $x$ to $y$. Since we do not want termination to be an issue we require that for all legal moves $x\rightarrow y$, $x>y$. For example $2\rightarrow 1$ or $2\rightarrow 0$ must be a move from N to P but, given the standard conditions for an invariant subtraction game, both $1\rightarrow 0$ and $3\rightarrow 1$ would represent $P\rightarrow P$, which is impossible. Therefore the set $A$ cannot be the set of P-positions for a standard subtraction game, as also follows from Theorem~\ref{t1}. However the situation can be remedied by the comply-number game from Section \ref{s3}.
\begin{Thm}\label{tt1}
Let $S = \{\{d,2d\}\mid d\in\N\}$. Then the set 
\begin{align*}
A &= \{0,1,3,4,9,10,12,13,\ldots \}\\ 
&= \left\{\sum_{i=0}^{\infty}\alpha_i3^i\mid \alpha_i\in \{0,1\}, \alpha_i=1 \text{ for at most finitely many } i\right\}
\end{align*}
is the set of all P-positions of the comply-number game $S$. It is the set of all N-positions for the comply-set game $S$.
\end{Thm}
\noindent{\bf Proof.} It is well-known that: the set $A$ can be obtained via an infinite greedy procedure which recursively includes the least nonnegative integer $z$ which does not satisfy $z+x=2y$ for any strictly smaller $x,y$ already in the set. (See \cite{KnLa04} for a proof of a generalization to $k$ term arithmetic progressions in case $k$ is prime.)

We will prove the first case. Then the second case follows from Theorem \ref{PN}. Thus, we have to show that a heap-size $x\not\in A$ is winning for the next player who declares the move set $s\in S$. By definition, this player can find a $d$ such that $\{x-d,x-2d\}\subset A$ for otherwise the greedy procedure would have included $x$ to the set $A$. By the comply-rule, the previous player, assigns one of these heap-sizes for the next player's move. In either case it is a P-position by induction. 

For the other direction, suppose that $x\in A$. Then, if the next player was able to announce a comply set, we have to demonstrate that the previous player can find a winning move. Since the greedy procedure included $x$ to the set $A$, the following statement is true: there is no $d>0$ such that $\{x-d, x-2d\}\subset A$. Suppose therefore that $\{x-d,x-2d\}$ is announced, but $x-d\not\in A$. Then the next player can be assigned, buy the previous player, to move to the heap-size $x-d$ (and otherwise to $x-2d\not\in A$). By induction such a heap-size is winning for the next player, that is the player not in turn to move from $x$. \hfill $\Box$\\

Thus our new construction solves the problem of finding a normal play impartial heap-game with $A$ (or $\M \setminus A$) as its complete set of P-positions. (In a sense we mimic the first board game in Section \ref{s2}.) It remains to investigate its nim-values (and analogously for the comply-set game). Also Mis\`{e}re versions (the last move loses) of our games may have interesting outcomes. But we will not pursuit these questions in this paper. 

An (invariant) \emph{game extension} of a comply game $S$ is a comply game $S\cup R$, for $R$ some set of finite subsets of $\N$.

\begin{Prob}\label{p2}
Is there any invariant game extension for the comply-number game with $S$ as in Theorem \ref{tt1}, with $A$ its set of P-positions, such that $$\#(A\cap \{0,\ldots , \frac{3^t-1}{2}\}) > 2^t$$ for any $t>0$ or for that matter such that $$\#(A\cap \{0,\ldots , n\}) \gg n^{\epsilon +\log2/\log 3}$$ for some $\epsilon > 0$?
\end{Prob}

Indeed, in \cite{KnLa04}, computer simulations proves that greedy gives the ``densest'' sets avoiding three-term arithmetic progressions for all subsets of $\{1,\ldots ,128\}$ in the following sense: the maximal number of integers within an interval of length $n\le 128$ which contains no three term arithmetic progressions is bounded above by $2^{\log_3(2n-1)}$. It is not very surprising that greedy is non-maximal for many $n$ since its distribution is very non-uniform. For example our computations show that there are denser sets than greedy for $n=85$, but nevertheless the number of numbers in this maximal set is $<2^{\log_3(2n-1)}$.

If there is such a game extension then its set of P-positions $A$ satisfies $x+z\ne 2y$ for all $x,y,z\in A$ where $x<y$. In other words the (heap-size) invariance of the set $S$ in the setting of Problem \ref{p2} implies that a set is void of three term arithmetic progressions. Our guess is that there is no such extension, but it would be interesting if this guess is wrong, because then it would be possible to construct ``denser'' sets than the greedy construction gives via invariant heap-games. Denser sets are known e.g. Behrend's famous \emph{finite} constructions where points from a $d$-dimensional hyper-sphere are projected to the natural numbers \cite{Be46}. Behrend's construction gives a very ``non-greedy'' set in one sense. It is much less dense than ``greedy'' for small sets, but grows much faster than greedy for larger sets. However the constructions are not comparable in another sense, since Behrend's approach is non-constructive. However, in \cite{Mo53}, Moser describes a constructive algorithm which gives finite sets avoiding three term arithmetic progressions with only a slighter less density then that of Behrend. His sequence begins $$100000, 1000100100, 1000400200, 1000900300, 1001600400, 1002500500, \ldots,$$ but for sufficiently large $n$ it contains at least $n^{1-c/\sqrt{\log n}}$ numbers in an interval of length $n$, for some constant $c$, hence the set obtained is much denser than that produced by greedy and only slightly less dense than Behrend's construction (which produces a smaller constant $c$). However it is easy to see that there can be no comply-subtraction game which has Moser's sequence as a set of P-positions, because each heap-size strictly between $100000$ and $1000100100$ must move to 100000 (which we treat as final position here). This holds also for the comply variation. But then for example 300100 is a legal move, which would connect two P-positions. Hence, in order to be able to use Moser's construction for an impartial game one needs some modification. Surely we can include the set $S$ as in Theorem \ref{tt1} as a subset of all rules and then adjoin \emph{non-invariant} rules which takes care of the remaining N-positions, but this is not necessarily a nice construction.

Anyway, it motivates a much more general definition of our games: for all heap sizes $x$, let the move sets $S_x$ depend on $x$. (Then $S=\{S_x\}$ is an invariant game if and only if $S_x=S_y$ for all heap sizes $x$ and $y$.) Otherwise, the game is defined, in analogy with the usual comply rules, by the family of sets of sets $\{S_x\}$: exchange $S$ for $S_x$ for each heap-size $x$. See also Section \ref{s51} for games on non-translation invariant arithmetic conditions.

We leave it as an open question to describe it or some nicer construction which emulates Moser's set as the set of P-positions. For example, does the list of non-invariant moves necessarily have to be infinite?

\begin{Prob}\label{p3}
Is it possible to obtain Moser's set (say with only finitely many numbers missing) with only a finite number of non-invariant moves?
\end{Prob}

In \cite{OdSt78} the authors explores small non-invariant disturbances to the greedy algorithm described in the first paragraph of the proof of Theorem \ref{tt1}: by shifting the first (non-zero) entry of the sequence but otherwise letting greedy generate the sequence as usual, they discuss whether the sets obtained are of the same order of magnitude or less dense than that of greedy, sometimes the patterns are regular, with fractal patterns similar to greedy, but for other initial disturbances they seem to be highly irregular. For example, for the case with the first entry $a_1 = 2$ the sequence begins $0,2,3,5,9,\ldots $. Their respective sets can be treated as sets of P-positions of games where, for each game, precisely one single non-invariant move is adjoined. They seem to suggest that any such disturbance gives a set of roughly the same upper asymptotic density as the non-disturbed greedy algorithm, although they only prove it for the regular expressions, that is whenever the first entry of the sequence is $3^i$ or $2\times3^i$ (here $i$ is a nonnegative integer). In our setting of games we can interpret for example the case $a_1=2$ as the invariant set $S$ together with the single non-invariant move from 1 to 0. Namely, then this sequence represents the complete set of P-positions for the slightly modified (almost invariant) comply-number game. In general one can adjoin any particular set of non-invariant moves to the set $S$ and thus obtain Moser's set as we have seen, which gives an answer to a question posed in that paper. This gives further motivation to Problem \ref{p3}. In \cite{ELRS99} the authors generalize the study in \cite{OdSt78} to a finite initial set of positive integers $A$, to obtain a ``Stanley sequence generated by $A$" and pose a number of interesting questions, for example (in our setting) the analogue of Problem \ref{p2} for finite non-invariant game extensions.

Remark: we note that Behrend's construction has recently been improved somewhat by Elkin, Wolf, Green, and generalized by O'Bryant but the gap to known upper bounds (Roth and generalizations/improvements by Szemer\'edi/Bourgain et al) is still enormous. It is generally believed that the upper asymptotic density of sets containing no three term arithmetic progressions is nearer Behrend's than (improvements of) Roth's result. 

\subsection{A certain game restriction}
On the other hand, one can remove move-sets from the set $S$ defined as in Theorem (\ref{tt1}) without affecting the status of the set $A$ as a complete set of P-position. We next show that there is a strict \emph{restriction} $T\subset S$ of this comply-number game such that its set of P-positions equals $A$. Indeed $T$ can be taken as generated by the set $A$ itself in the following sense.

\begin{Thm}\label{ttt}
Let $A=\{0,1,3,4,\ldots \}$ be as in Theorem \ref{tt1}. Then $A$ is the complete set of P-positions of the comply-number game $T=\{\{d,2d\}\mid d\in A\}$. That is, if we apply the greedy algorithm to avoid precisely three term arithmetic progressions of discrepancies $d\in A$, then we produce the set $A$. On the other hand, let $U$ denote some strict restriction of $T$. Then the set of P-positions of $U$'s comply-number game contains a three-term arithmetic progression.
\end{Thm}

\noindent{\bf Proof.} Write $x\not\in A$ in base three expansion. Then $x=\sum_{i\in \M} x_i3^i$ for appropriate choices of $x_i\in \{0,1,2\}$. Then let $d= \sum d_i3^i$, where $d_i=1$ if $x_i=2$ and $d_i=0$ otherwise. This gives $d>0$ and the next player can find a subset of $A$ to move to. Also there can be no move between positions in $A$ since $T$ is a restriction of $S = \{\{d,2d\}\mid d\in \N\}$.

For the second part, suppose that $d = \sum_{i\in \M} d_i3^i$, $d_i\in\{0,1\}$, is the least discrepancy such that $\{d,2d\}\in T$ but $\{d,2d\}\not \in U$. Let $x = \sum_{i\in \M} x_i3^i$ be such that $x_i = 2d_i$ for all $i$. Then for all $c = \sum_{i\in \M} c_i3^i$ there is a least $i$ such that ($x_i=2$ and $c_i=0$) or ($x_i=0$ and $c_i=1$). If the first case holds then (since the second case is ruled out) there is a digit $x_i-c_i=2$. If the second case holds, there will be a carry from some $x_j=2$, with $j>i$ minimal. Hence $x_i-c_i=2$ also in this case. Hence there is no move-set which takes $x$ to a subset of $A$. But, by minimality of $d$, $A\cap [0,x)$ is identical to the set of numbers $<x$ produced greedily by $U$. This produces the arithmetic progression $x-2d,x-d,x$ of which, by definition of $d$ and the previous sentence, each number is in this latter set. \hfill $\Box$\\

\begin{Prob}
Is it possible to find any set $D$ with upper asymptotic density $o(n^{\log2/\log 3})$ such that the set of P-positions of the comply-number game defined by  $\{\{d,2d\}\mid d\in D\}$ does not contain any three-term arithmetic progressions (of any discrepancy)?
\end{Prob}

\subsection{Another $\star$-operator and self duality}
Let $D\subset \N$ and let the move-sets of the comply-number game $S^\star$ be defined by the set of P-positions $A$ of $S=\{\{d,2d\}\mid d\in D\}$ in the following way $S^\star=\{\{a,2a\}\mid a\in A\}$. By the result in Theorem \ref{ttt} we get 
\begin{Cor}\label{cor}
Let $D=A$ with $A$ as in Theorem \ref{tt1}. Then $S^\star=S$ and hence, for all $k$, $S^{k\star}=S$. 
\end{Cor}
Thus we get self-duality for this $\star$-operator a result to be compared with the results in \cite{LHF11, La12} for subtraction games (on several heaps) $\mathcal{M}$ where self-duality never holds (but sometimes $\mathcal{M}=\mathcal{M}^{\star \star}$). For what sets $D$ do we get $S=S^\star$? For what arithmetic conditions and given an appropriate definition of ``discrepancy'' do we have analogues to Corollary \ref{cor}?

\subsection{Other interpretations}\label{s43} 
The three-term arithmetic progression game in Section \ref{s2} is in a sense only a reformulation of the comply-number rules, played as a board game. Is there also an invariant heap-game without comply restrictions that solves the set $A$ (and similar sets) from Section \ref{s2}? The following game is played on three heaps of tokens, represented by a triple of heap-sizes $(x,y,z)$ with $x < y $ and $x+z=2y$: the next player removes the largest heap and one of the smallest heaps and let the remaining heap become the largest in the next position, say for example the heap with $y$ tokens. In addition he chooses a positive integer $d$ and presents the position $(y-2d,y-d,y)$ for the next player, provided $y-2d\ge 0$. In this way the number of tokens will strictly decrease for each move. The final position, which is in P, will be $(0,1,2)$. We get the following result:
\begin{Thm}\label{ttt1}
The position $(x,y,z)$ is in N if $z\in A=\{0,1,3,4,\ldots \}$ as defined in Theorem \ref{tt1}. Otherwise, if $z\not\in A$, then $(x,y,z)$ (where $x< y<z$) is in P if and only if $\{x,y\}\subset A$.
\end{Thm}
The proof is similar to that of Theorem \ref{tt1}, and uses again that the set $A$ is the greedy construction of a set which does not contain three term arithmetic progressions, so we omit it. Therefore, by abuse of notation, we can simplify the statement of Theorem \ref{ttt1} and say that $(x,y,z)$ is in N if and only if $z\in A$, because if $z\not\in A$ then, by greedy, it is always possible to choose $x,y$ appropriately so that $\{x,y\}\subset A$.

For other arithmetic constraints we can use the number of variables in the constraint to represent the number of heaps in the game and proceed in analogy to the case for three term arithmetic progressions.

\section{One heap games in general}\label{s5}

We will have use for a general definition of our conditions. 
A \emph{linear form} is here an expression on finitely many variables with integer coefficients $f(x_1,\ldots ,x_k) = \alpha_kx_k + \ldots + \alpha_1x_1 + \alpha_0$. Then $f$ is \emph{translation invariant} if and only if $\alpha_1+\ldots +\alpha_k=0$. Suppose that we have a finite or countable family of linear forms $F(x_1,\ldots ,x_k)=\{f_i(x_1,\ldots ,x_k)\}$. For each $i$, let $e_i =e_i(x_1,\ldots ,x_k) $ be a boolean variable which assesses whether there exists a non-trivial solution to $f_i(x_1,\ldots ,x_k) = 0$. 

A \emph{trivial solution} $(x_1,\ldots ,x_k)$ satisfies $\alpha_0 = 0$ and $\sum_{i\in I} \alpha_i = 0$, for all sets of indices $I\subset \{1,\ldots ,k\}$ such that $x_i=x_j$, for all $i,j\in I$. Thus trivial solutions can only exist for translation invariant conditions (see also Section~\ref{s51}). 

Given a family of linear forms $F$, let $AC=AC(F(x_1,\ldots ,x_k))$ be the (possibly uncountable) family of (possibly infinite) expressions consisting of the $e_i$s, the connectives AND and OR and well-formed parenthesizes. Then a set $X$ of integers \emph{avoids} $ac\in AC$ if $ac$ is false for all $k$-tuples in $X$. For the purpose of this paper, we assume that ac is decidable. We will return to these general definitions in Section \ref{s6} and illustrate their use in a few examples (while the main reason for going this general is the setting of Theorem \ref{thm42}). 

Similar greedy approaches can be used in combination with other arithmetic constraints, producing sets that avoid non-trivial solutions to systems of linear equations in finitely many variables e.g. \cite{Ru93,Ru95} (system of \emph{linear forms}). 

For examples of \emph{translation invariant} systems, for a given integer $k\ge 2$, we have: 
\begin{itemize}
\item the \emph{Sidon condition} which avoids repetitions of $x_1 + \ldots + x_k$ ($2k$ variables). A trivial solution here is simply a permutation of the indices.
\item the \emph{arithmetic mean condition} avoiding solutions to $x_1 + \ldots + x_{k-1} = (k-1)x_k$ ($k$ variables). Here triviality means that all entries are the same.
\item the condition which avoids solutions to non-trivial ($d>0$) $k$-term arithmetic progressions $x,x+d,\ldots ,x+(k-1)d$, which means avoiding simultaneous non-trivial solutions to a system of $k-2$ equations of the form $x_i+x_{i+2}=x_{i+1}$, $i\in \{1,\ldots , k-2\}$ (here one equation is non-trivial if and only if all are). 
\end{itemize}
There are base $k$ constructions corresponding to the latter two greedy formulations, but none are known for the Sidon conditions. For arithmetic mean avoidance in $k$ variables, it is not hard to check that the greedy algorithm produces a set that uses only the digits $0,1$ in base $k$; for extensive generalizations of this condition with the same base $k$ construction, see \cite{Ts11}; for $k$-term arithmetic progressions for $k>3$ a prime, never use digit $k-1$ except for the least digit where we never use 0, see also \cite{KnLa04}. Of course, generalizations of  Theorem \ref{ttt} hold for these constructions. We note that there is a diversity of so-called ``nonaveraging" conditions in \cite{Ts11}---we will not go into more detail here since we only recently learned about this beautiful paper---that satisfy one and the same base $k$ construction. This implies that there are in general many invariant comply games for one and the same such set of P-positions. This appears to be a fruitful area for future studies. 

Any set avoiding $k$-term arithmetic progressions must have upper asymptotic density zero as Szemer\'edi demonstrated in \cite{Sz75}. An analogous result holds for any arithmetic condition if and only if it is translation invariant \cite{Ru93, Ru95} (see also \cite{KnLa04} for a discussion).

\subsection{Some non-translation invariant conditions}\label{s51}
For a non-translation invariant condition in one equation one can take for example 
\begin{align}\label{kyxz}
ky=x+z 
\end{align}
for $0<k\ne 2$. In \cite{BHKLS05} very explicit finite constructions are studied for sets of "maximal" density avoiding (\ref{kyxz}). We do not know whether greedy constructions have been systematically studied for such conditions, except for some simple cases. Sets avoiding equations in two variables has been studied in \cite{KnLa04} where a greedy algorithm in fact is shown to give sets of maximal cardinalities (but these do not appear to be very interesting as games) where also a certain ``semi-greedy" algorithm is proved to produce least dens so-called saturated sets. We note that the greedy algorithm produces a maximal set also for the case $k=1$ in (\ref{kyxz}). Namely 1 will be included, but 2 will not since $2=1+1$, then 3 is OK since 2 wasn't included and so on, which gives the set of odd numbers, $\{1,3,5\ldots \}$. Note that we start the greedy algorithm at ``1'' rather than ``0'' here, to avoid trivialities. It would possibly be more convenient to always work in $\N$ rather than $\M$, but in the setting of heap games it is customary to identify the empty heap with ``0'' and our main interests in this paper are translation invariant conditions. 

The sets defining non-invariant games for for a family $F$ of non-translation invariant linear forms are particularly simple: $S_x=S_x(ac) = \{\{x_i\}_x\mid x_i<x, 0<i\le k, e_i \text{ holds }\}$, where $x = x_j$ for some $j\ne i$ and $\{x_i\}_x$ denotes the unordered $(k-1)$-tuple resulting from omitting $x$. For example, for $F=\{x+y=3z\}$, we can take $S_x=\{\{y,z\}\mid x+y=3z\}$. The rules appear attractive although non-invariant. It would be interesting to study  $ac=\{x+z=2y \text{ OR }x+z=3y\}$ and similar conditions, that is ``disturbing" greedy avoiding three term arithmetic progression by a non-translation invariant condition. Note here that the variables in the two equations are independent, but for a corresponding ac with ``OR" exchanged for "AND" the situation is more diverse. Namely there are several conditions possible from two equations, for example $ac=\{x+z=2y \text{ AND }x+z=3y\}$ does not have any solution on $\N$, whereas $ac=\{x+z=2y \text{ AND }x+z=3w\}$ has infinitely many. Hence the number of variables for an $ac$ depends on the particular connectives used: $\#(X \text{ OR } Y)=\#\max \{\#X,\#Y\}$ whereas $\#(X \text{ AND } Y) = \#\{X,Y\}$. When are such greedy definitions regular in the sense of \cite{OdSt78}?

We believe that the following problem has a positive answer.

\begin{Prob}
Has each greedily produced set of numbers given a non-translation invariant avoidance criterion positive upper asymptotic density? 
\end{Prob}

The question is motivated by the observation that each set of maximal cardinality avoiding a non-translation invariant condition has a positive asymptotic density. A related question is wether such conditions eventually become periodic.

\begin{Prob}
Does each greedily produced set of numbers given a non-translation invariant avoidance criterion eventually become periodic? 
\end{Prob}

Take, for example $k=3$, in (\ref{kyxz}). The greedily produced set begins $1,3,4,7,10,12,13,15,16,19,22,25,\ldots$. It is not obvious from the initial elements whether this sequence will eventually become periodic. However one can prove that all numbers congruent to 1 modulo 3 belong to the sequence, but none congruent to 2, and also that infinitely many are divisible by 3. As a second observation, it is not possible to define an invariant comply game with this sequence representing the set of P-positions (terminal position is heap-size $1$). Namely it is forced that $\{1\}\in S$, but this connects 4 and 3, which is impossible if they were both in P. Note, however that this is possible for the set $\{1,3,5,\ldots \}$ produced by the non-invariant condition $x = y + z$. Thus, sometimes a non-invariant condition permits invariant move sets, but in general this is not possible. 

\begin{Prob}
Classify the non-translation invariant arithmetic conditions that permit invariant game rules. If a particular condition does not permit invariant rules, is it possible to restrict the non-invariant part to a finite set?
\end{Prob}
 
\section{Comply games on several heaps}\label{s6}
One can check that results in analogy to those in Section \ref{s3} will still hold if we let the symbols represent $d$-dimensional vectors (exchange $\M$ for ${\M}^d$) and where inequalities are interpreted as usual for partially ordered sets that is $a<b$ with $a,b\in \M^d$ means that $a_i\le b_i$ for all $i\in\{1,\ldots ,d\}$ and with strict inequality for at least one component. Also $a-b=(a_1-b_1,\ldots ,a_d-b_d)$, where, for all $i$, $a_i-b_i\ge 0$, but in Section \ref{s3} we kept the terminology to a minimum as not to obscure the ideas, e.g. \cite{La12}. Each subtraction set in $S$ is a set of vectors of nonnegative integers, at least one of them strictly greater than zero. That is, given a position $(x_1,\ldots ,x_k)=x\in \M^k$ and a subtraction set $S\in \M^k$, the \emph{move} $x-s$ is legal if and only if $(s_1,\ldots ,s_k)=s\in S$ and for each $i$ $0\le x_i-s_i\le x_i$ with at least one $x_i-s_i<x_i$. 

How about the results in Section \ref{s4}? It is not immediately clear how one would generalize the greedy rule for $d>2$, but the case $d=2$ is studied in \cite{He04, KnLa04}. Even here there are several choices. A function $\pi:\M\rightarrow\M$ contains a three term arithmetic progression if there exist $x,y,z\in \M$ such that $x+z=2y$ with $x<y$ and $\pi(x)+\pi(z) = 2\pi(y)$; it \emph{avoids three term arithmetic progressions} if there is no such triple $(x,y,z)$, \cite{Si88}. In \cite{He04} a greedy injection $\pi_g:\M\rightarrow\M$ avoiding three term arithmetic progressions is defined, by letting $\pi_g(n)$ be the least positive integer such that 
\begin{align}\label{mn}
\pi_g(m) \ne \pi_g(n) 
\end{align}
for all $0\le m<n$ and such that $\pi_g(n)$ does not form an arithmetic progression with any previous entries. Moreover, it is demonstrated that $\pi_g$ is a \emph{permutation} with, for all $n$, $3/8 \le \pi(n)/n < 3/2$, but otherwise its behavior is not yet well understood, for example we do not yet know whether or not $\pi_g(n)/n$ converges (to $1$). (An intermediate result following from Szemer\'edis theorem is that $\pi_g(n)=n+\Omega(1)$, \cite{KnLa04}.) If one removes the requirement (\ref{mn}) then we get a function with equally interesting properties \cite{KnLa04}. An open question is if the analogous greedy injection for the Sidon condition is a permutation. To our awareness this question has not been settled yet. This motivates our main result in this section. 

We will now define some variations of $\pi_g$. Let the greedy choice disregard any strictly decreasing arithmetic progressions, two choices here: with or without (\ref{mn}). We will next show that the injection $\pi_g$ (including (\ref{mn})) thus obtained is a permutation, more precisely an \emph{involution}, that is satisfying $\pi_{g}(n)=\pi_{g}^{-1}(n)$ for all $n$. We conjecture that for this case it holds that $\pi_g(n)/n\rightarrow 1$, motivated by that the arithmetic constraint is weaker for this case. 

Before we prove that this variant of $\pi_g$ is an involution, let us first generalize the concept of sets avoiding arithmetic conditions to functions avoiding arithmetic conditions. (We wish to find a more general class for which $\pi_g$ is an involution.)

Given a family of linear forms $F=F(x_1,\ldots ,x_k)$, let $ac\in \text{AC}(F)$. Then we say that a function $\pi:\M\rightarrow \M$ (or $\pi: \N\rightarrow \N$ in case of non-translation invariant conditions) avoids $ac$ if, for all $(x_1,\ldots ,x_k)$ such that $ac(x_1,\ldots x_k)$ holds, $ac(\pi(x_1),\ldots \pi(x_k))$ does not hold. 

As we have already remarked, the logical expressions may not be decidable since we include the possibility of infinitely many linear forms. Hence it is not immediately clear whether one can decide if, given an $ac$, there exists a function which avoids $ac$. 
But, since only finitely many equations need to be tested for a greedy definition (there are only finitely many coefficients less than a given $n$ for each of $\pi_g$'s tests) this question is settled. 

\begin{Thm}
Given any $ac$ there is a function that avoids $ac$. One such function is $\pi_g$.
\end{Thm}

This result carries over to any \emph{relaxed} avoidance (or equivalently restricted arithmetic condition). Some relaxation are motivated by nice game rules to come. A function $\pi: \M\rightarrow \M$ avoids $ac$ weakly if, for all $x_1\le \ldots \le x_k$ such that $ac(x_{1},\ldots x_{k})$ holds, $ac(\pi(x_{1}),\ldots \pi(x_{k}))$ does not hold whenever $\pi(x_1)\le \ldots \le \pi(x_k)$. Thus, for this case, we only avoid a progression if $\pi$ is \emph{order preserving} on this progression. Another restriction of an arithmetic condition is that whenever $x_j = \max_i\{x_i\}$ then $\pi(x_j) = \max_i\{\pi(x_i)\}$. Hence this class, called here \emph{$\max ac$}, gives a stronger avoidance than order preserving, but weaker than for unrestricted $ac$. Thus, if a function $\pi$ avoids $ac$ then it avoids $ac$ weakly and $\max ac$. If it avoids $\max ac$ then it avoids $ac$ weakly.
Hence let us first reformulate our first example on the greedy injection avoiding 3-term arithmetic progressions in the new general context of arithmetic conditions  relaxed by the $\max ac$ condition, which will be our primary interest. We define $\pi_g: \M\rightarrow \M$ recursively. Given $\pi_g(0)=0$ and $\pi_g(m)$ for all $0\le m<n$, let $\pi_g(n)$ be the least positive integer such that $\pi_g(m) \ne \pi_g(n)$ for all $m$ and such that for all 
\begin{align}\label{xn1}
\{x_1,\ldots ,x_{k-1}\}\subset \{1,\ldots ,n-1\} 
\end{align}
such that 
\begin{align}\label{x}
ac(x_{1},\ldots, x_{j},n,x_{{j+1}},\ldots, x_{{k-1}}) 
\end{align}
holds then 
\begin{align}\label{pi}
ac(\pi_g(x_{1}),\ldots,\pi(x_{{j}}),\pi_g(n),\pi_g(x_{{j+1}}),\ldots, \pi_g(x_{{k-1}})) 
\end{align}
does not hold whenever 
\begin{align}\label{<}
\pi_g(x_i) < \pi_g(n) 
\end{align}
for all $i\in \{1,\ldots ,k-1\}$. By this definition it is clear that $\pi_g$ is an injection that avoids $\max ac$. The main result in this section demonstrates that, in fact, $\pi_g$ is an involution. 

\begin{Ex}\label{e1}
The rules of the classical game of Nim (on $k\in \N$ heaps of tokens) are as follows: a player in turn removes any number of tokens, at least one and at most a whole heap form precisely one of the heaps. The empty arithmetic condition gives Nim's P-positions playing on two heaps of tokens, that is for all $n$, $n=\pi_g(n)$. We will include the Nim type moves to all our games since, as we will see, they correspond to $\pi_g$ being a function satisfying (\ref{<}) (that is an injection). 
\end{Ex}

We have chosen the following rules of game. Playing from position $(x,y)$, the \emph{max-ac comply rules} are defined by the set $S=S(ac)$ which includes all Nim-type moves and otherwise ranges over those unordered $(k-1)$-tuples $((x_i,y_i))$ that together with $(x,y)$ satisfy $ac$ and where each entry of the tuple satisfies $x_i < x$ and $y_i < y$. For the order preserving condition the sets in $S(ac)$ are further restricted. Namely, for all $i\ne j$ such that $x_i<x_j$ we require $y_i<y_j$. For the unrestricted $ac$ (where only trivial solutions are removed) we only require that $x_i<x$ for all $i$. As we have remarked it is not always clear whether a given injection $\pi_g$ is a permutation. We get the following result.

\begin{Thm}
A position of the comply-number game $S(ac)$ is in P if and only if it is of the form $(x,\pi_g(x))$  or $(\pi_g(x),x)$. 
\end{Thm}
\noindent{\bf Proof.} Since $\pi_g$ is an injection, a Nim-type move can never connect two positions of the given form, by definition of greedy and (\ref{<}). If $(x,y)$ is of the form $y>\pi_g(x)$ then there is not necessarily a set $s\in S$ with all entries of the given forms. But then a Nim type move suffices to find such a position, namely $(x,\pi_g(x))$. Otherwise, if $y<\pi_g(x)$ it suffices to find a set $s$ with each position of the given form. But such a set exists, unless there is an $n<x$ such that $y=\pi_g(n)$ in which case a Nim-type move suffices, by the definition of $\pi_g$, for otherwise greedy would have chosen $y$. This gives the other direction as well, since, by definition, there can be no such set whenever $y=\pi_g(n)$. 
\hfill $\Box$\\

\begin{figure}[ht!]
\begin{center}
\vspace{0.5 cm}
{\includegraphics[width=0.49\textwidth]{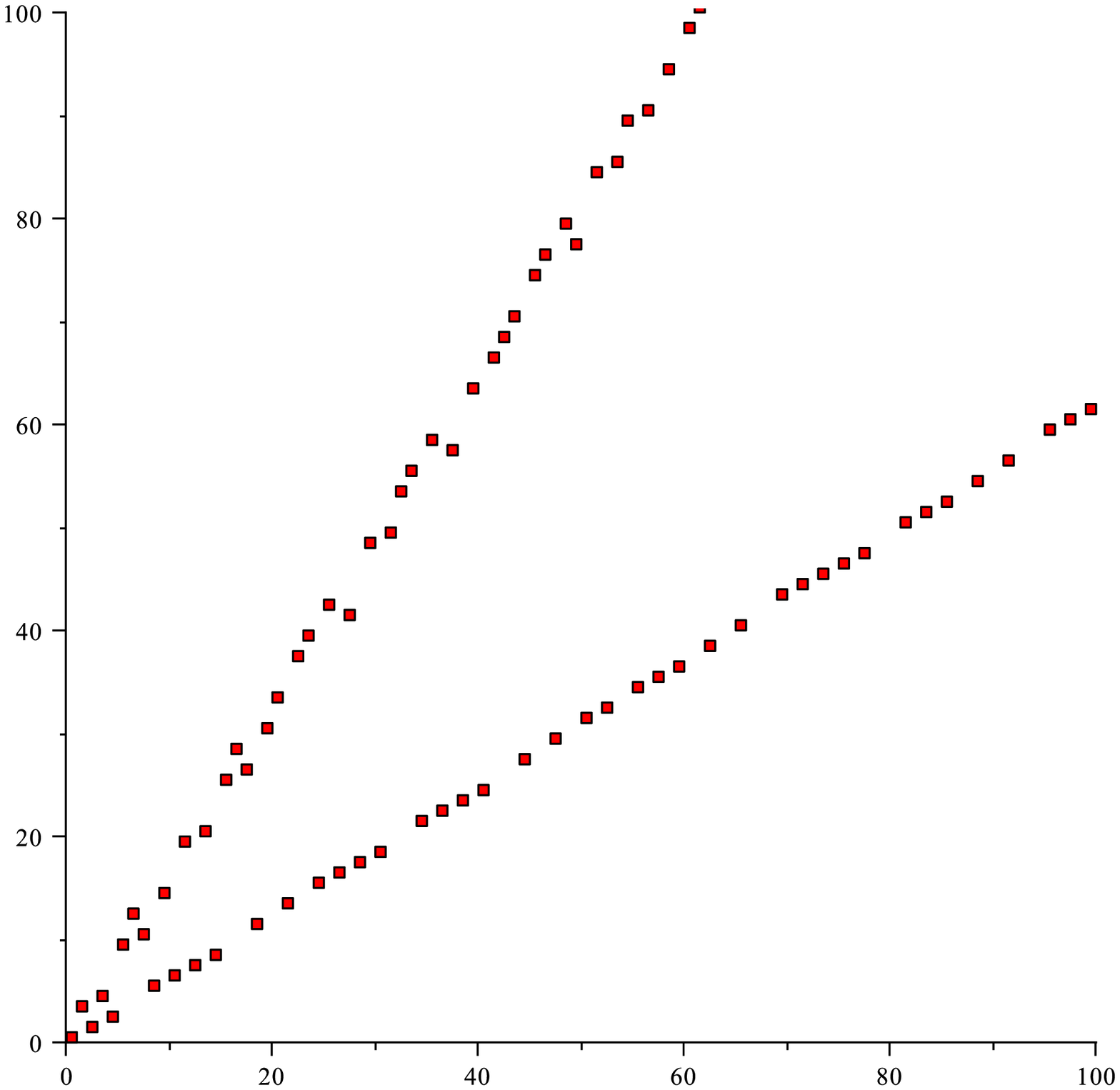}}\hspace{0.1 cm}{\includegraphics[width=0.49\textwidth]{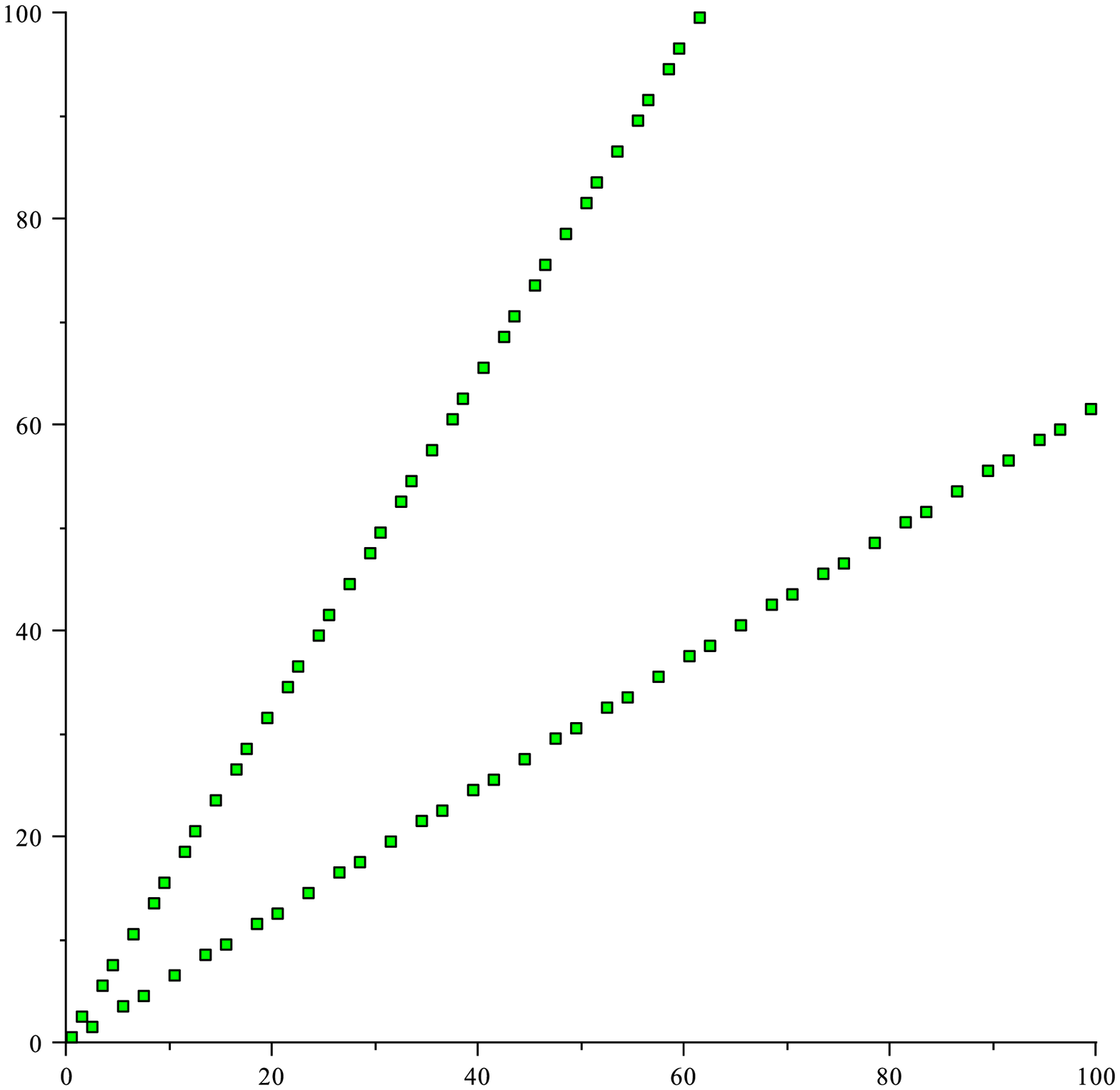}}
\end{center}\caption{Asymmetric and classical Wythoff Nim with initial P-positions $(0,0),(1,3),(2,1),(3,4),(4,2),(5,9),(6,12)$ and $(0,0),(1,2),(2,1),(3,5),(4,7),(5,3),(6,10)$ respectively.}\label{W}
\end{figure}

\begin{figure}[ht!]
\begin{center}
\vspace{0.5 cm}
{\includegraphics[width=0.49\textwidth]{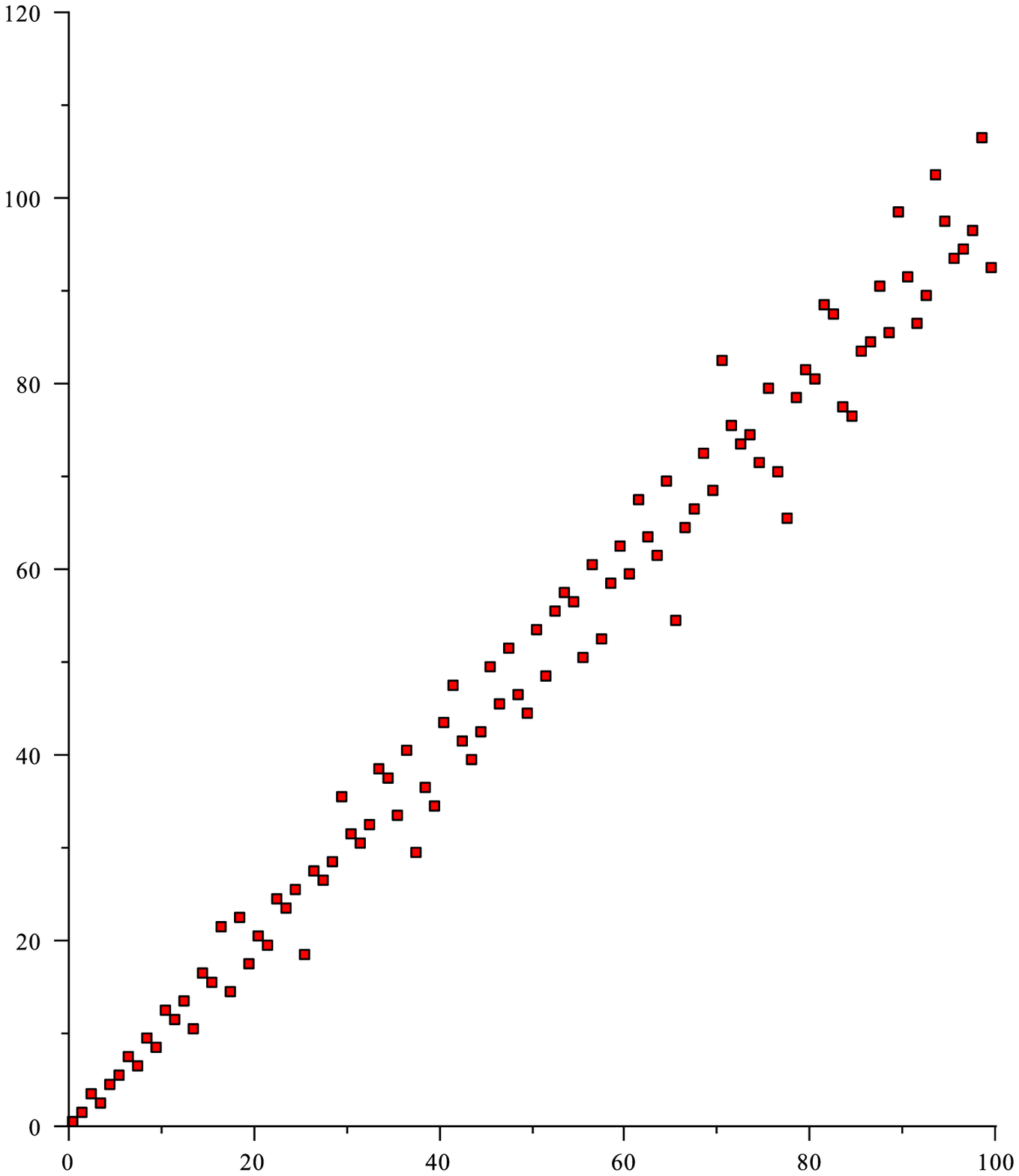}}\hspace{0.1 cm}{\includegraphics[width=0.49\textwidth]{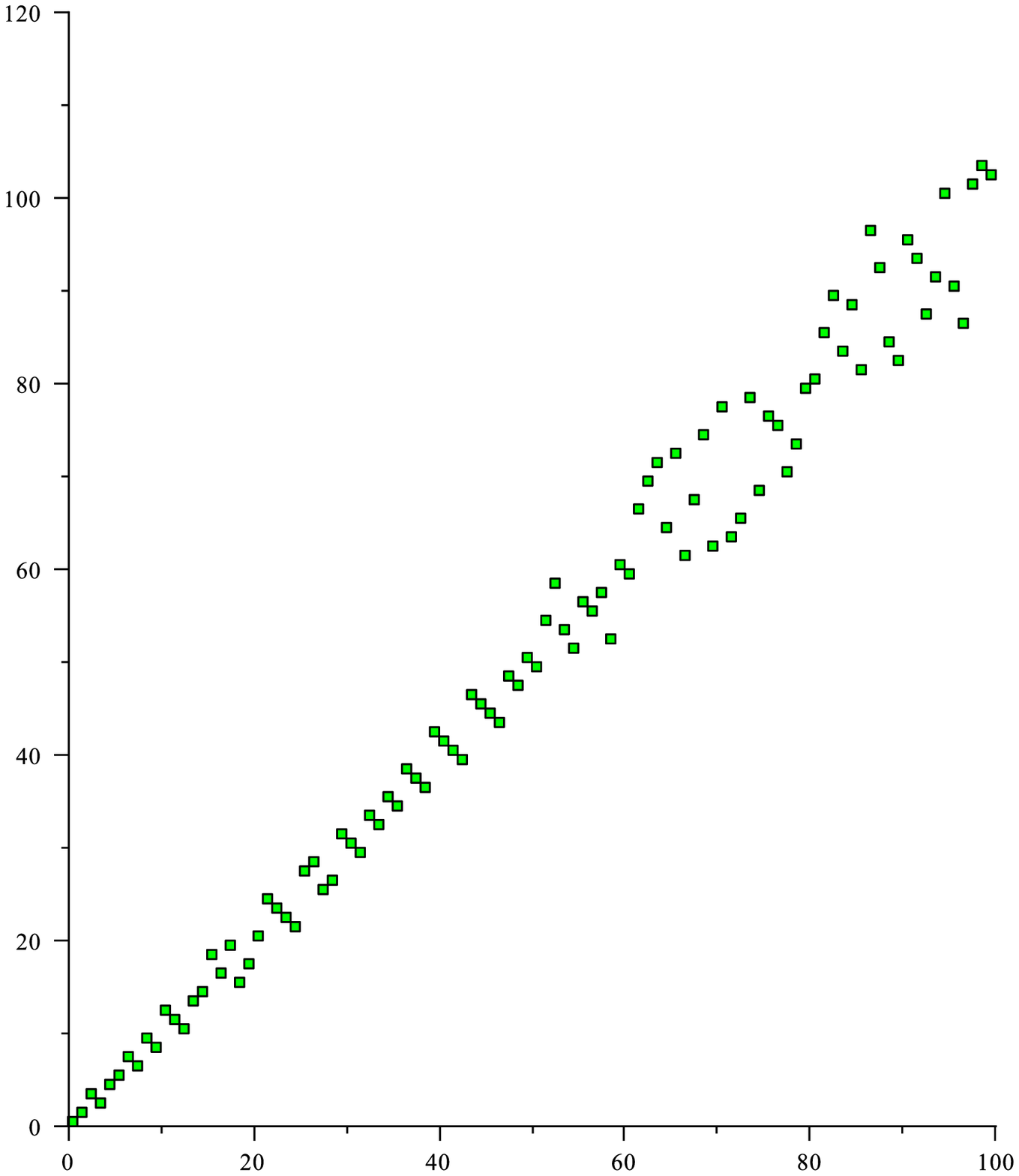}}
\end{center}\caption{Asymmetric (unrestricted avoidance) and symmetric (max $ac$) 3-term-greedy. }\label{AP}
\end{figure}

\begin{figure}[ht!]
\begin{center}
\vspace{0.5 cm}
{\includegraphics[width=0.25\textwidth]{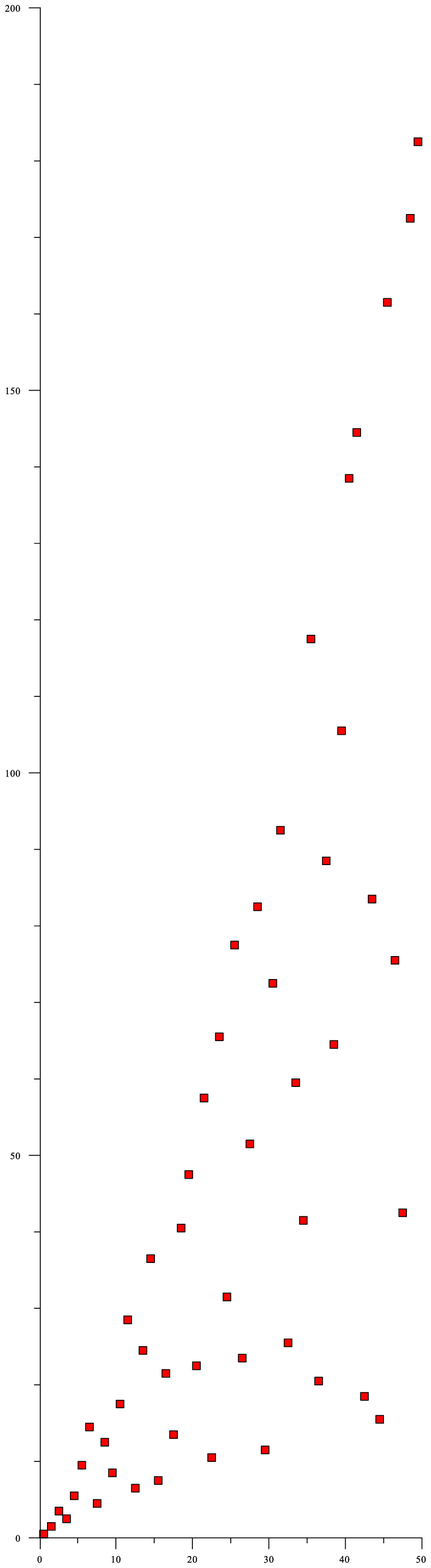}}\hspace{0.5 cm}{\includegraphics[width=0.25\textwidth]{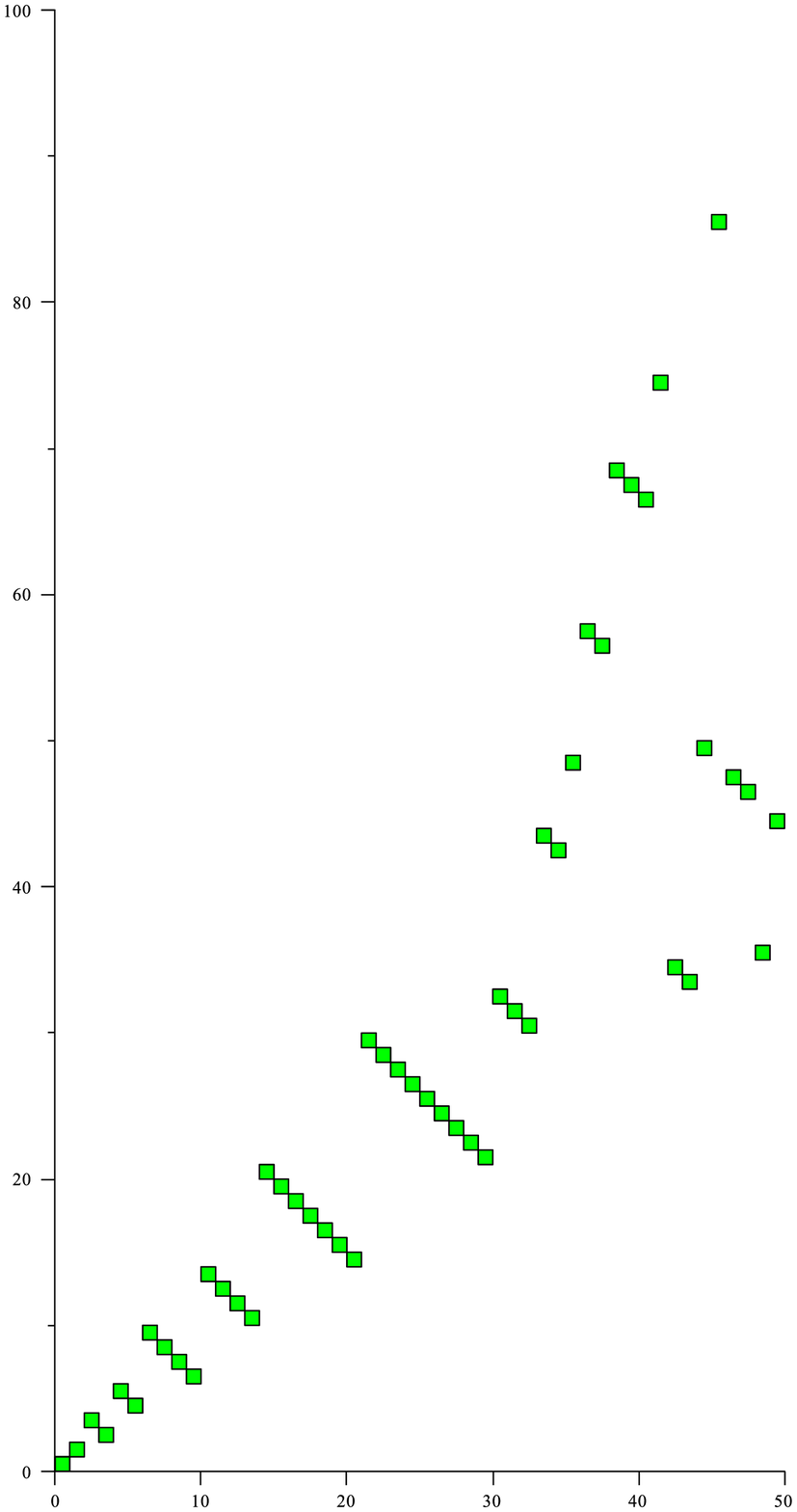}}\hspace{0.5 cm}{\includegraphics[width=0.25\textwidth]{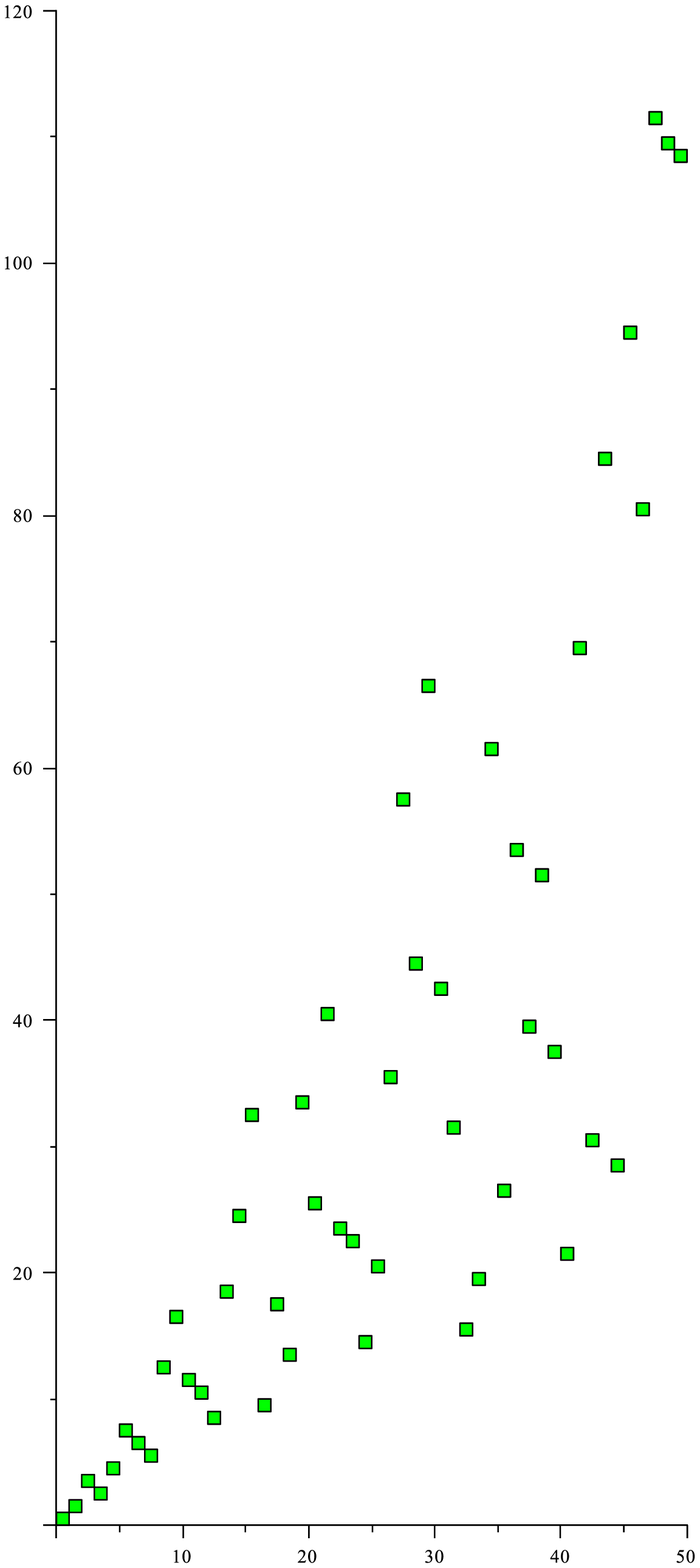}}
\end{center}\caption{Asymmetric (unrestricted) and symmetric (max $ac$ and order preserving $ac$ respectively) Sidon-greedy.}\label{Sid}
\end{figure}

\begin{figure}[ht!]
\begin{center}
\vspace{0.5 cm}
{\includegraphics[width=0.49\textwidth]{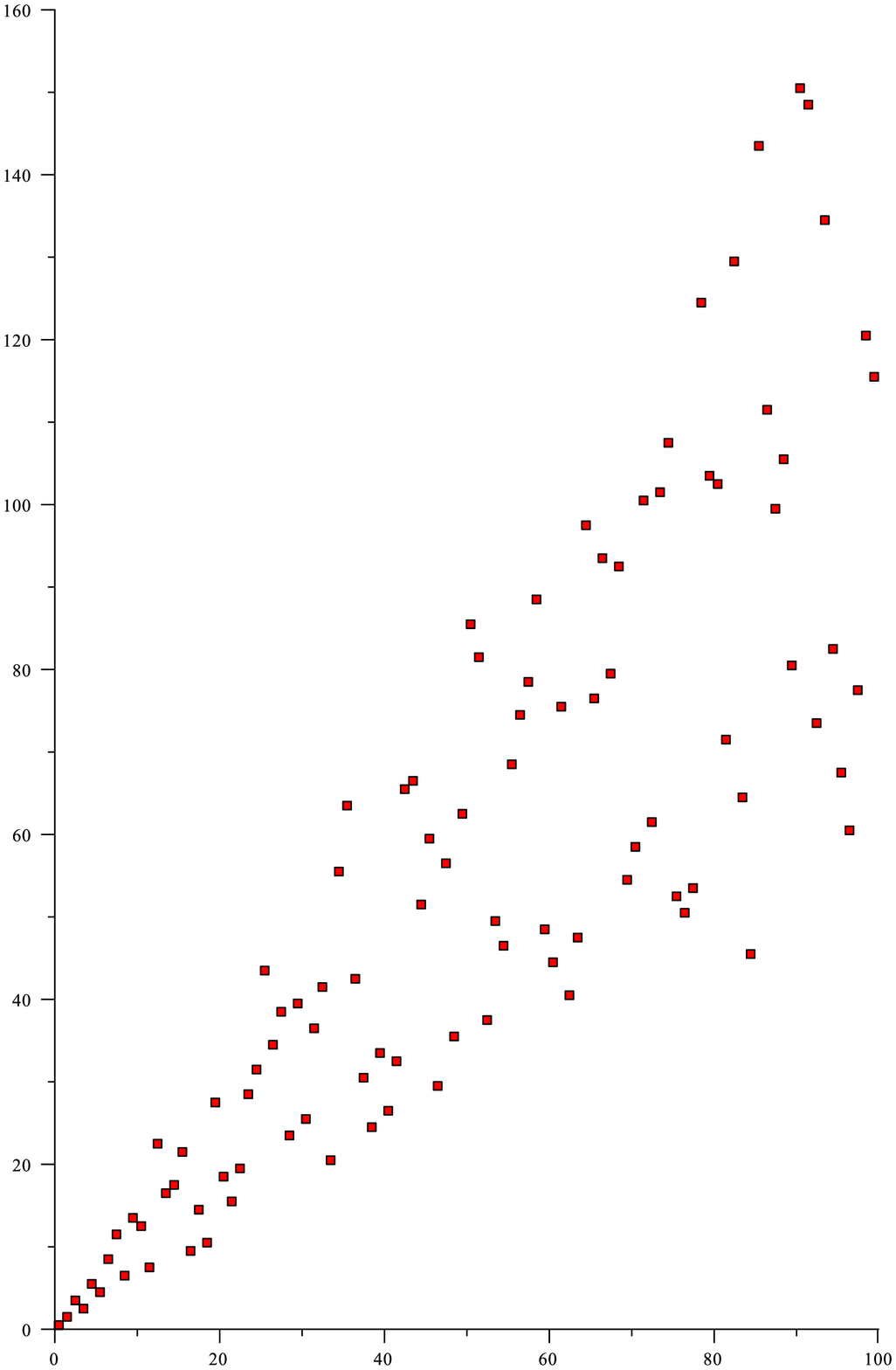}}\hspace{0.1 cm}{\includegraphics[width=0.49\textwidth]{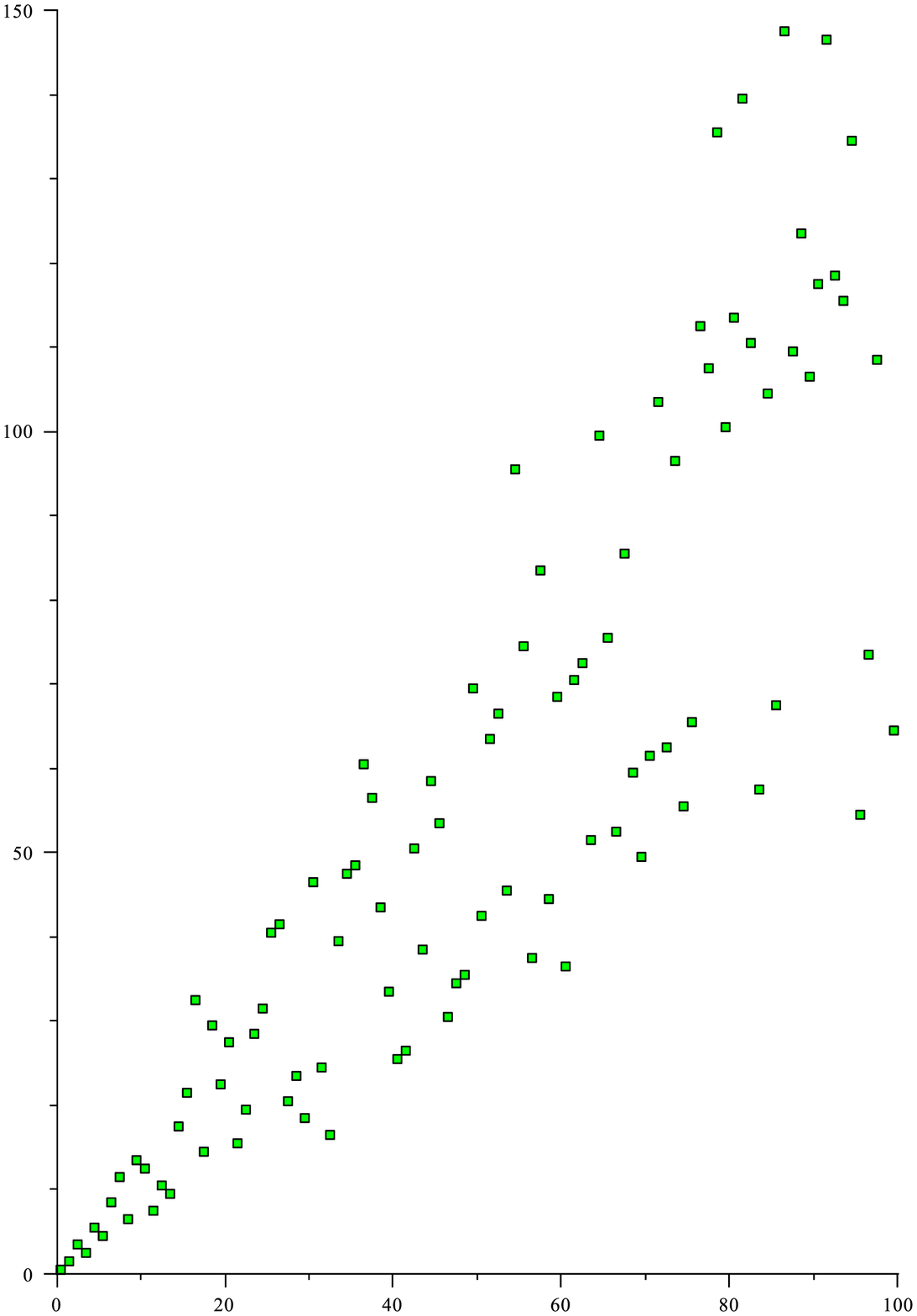}}
\end{center}\caption{Asymmetric and symmetric line-greedy.}\label{Line}
\end{figure}

\begin{figure}[ht!]
\begin{center}
\vspace{0.5 cm}
{\includegraphics[width=0.24\textwidth]{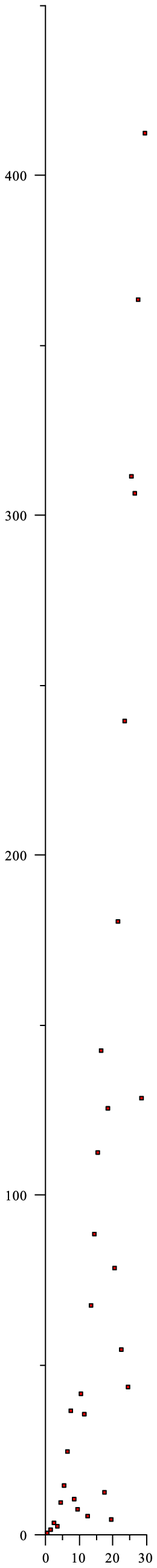}}\hspace{0.5 cm}{\includegraphics[width=0.24\textwidth]{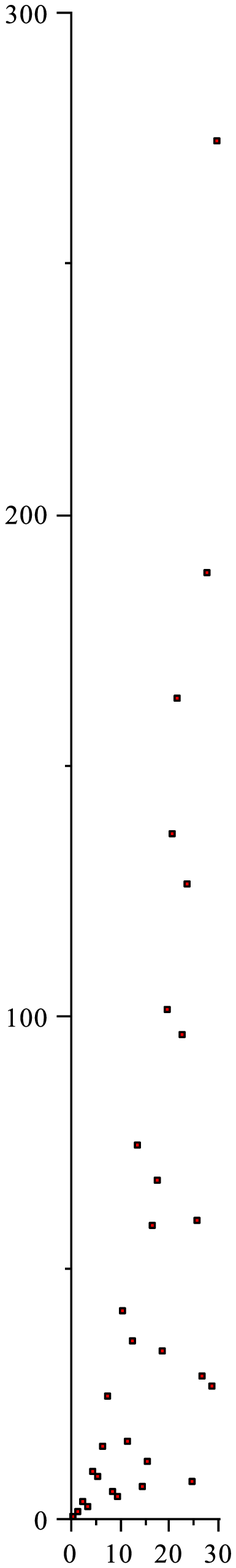}}\hspace{0.5 cm}{\includegraphics[width=0.24\textwidth]{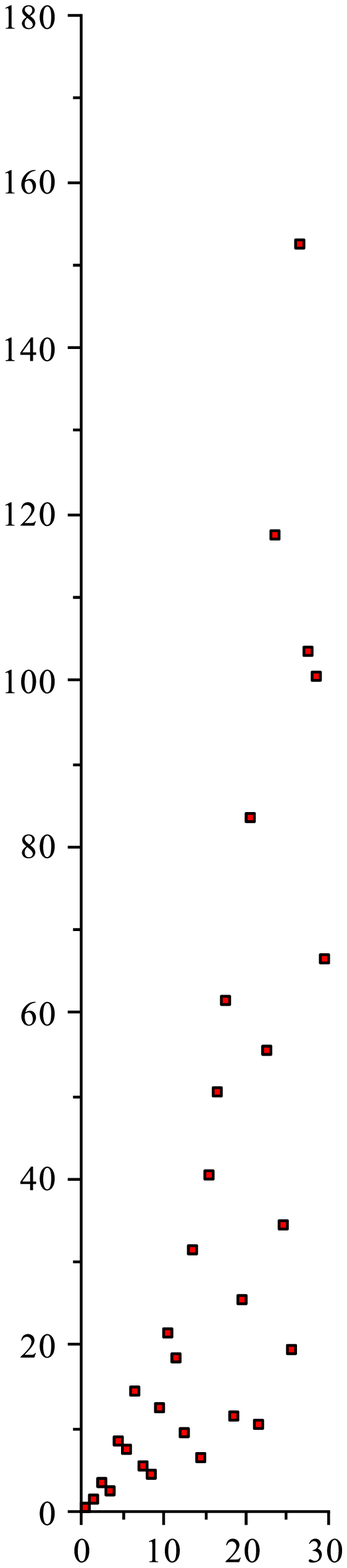}}
\end{center}\caption{Asymmetric (unrestricted) and symmetric (max $ac$ and order preserving $ac$ respectively) parallel-greedy.}\label{para}
\end{figure}

Let us view the following examples: \emph{Wythoff Nim, Sidon-greedy, k-term-greedy, line-greedy}, also illustrated in Figures \ref{W} to \ref{para}. For line-greedy, see also Section \ref{s2}.
\newpage
\begin{Ex}
To generate the P-positions of \emph{Wythoff Nim}, we can take $f_i=x_2-x_1-i$, for all $i\in \N_0$ and $ac = e_1 \text{ OR } e_2 \text{ OR } e_3 \text{ OR }\ldots$. Therefore, for each $n=x_2$, it suffices to test $x_1$ for $1,\ldots , n-1$ (although there exist faster algorithms). Note that the P-positions will depend on the choice of restriction of $ac$. For the classical P-positions we chose for example max $ac$. The unrestricted variant (we call it \emph{asymmetric Wythoff Nim}) does not give an involution, but we believe one can prove that indeed, it is a permutation and therefore corresponds to the games we have defined. It is an interesting open problem to investigate the P-positions of this game. We can also think of it as: those bishop type moves that increase the second coordinate are also legal (in addition to the usual rules). 
\end{Ex}
\clearpage
\begin{Ex}\label{e3}
For $k$-term-greedy we can take $f_i=x_i+x_{i+2}-2x_{i+1}$ for $i\in \{1, \ldots , k-2\}$ and $ac = e_1 \text{ AND } e_2 \text{ AND }\ldots \text{ AND } e_{k-2}$. For $k=2$, both unrestricted and max $ac$ has initial P-positions $$(0,0),(1,1),(2,3),(3,2),(4,4),(5,5),(6,7),(7,6),(8,9),(9,8),(10,12),(11,11).$$ The next entry is the first that differs, $(12,10)$ and $(12,13)$ respectively. See also \cite{KnLa04} for an extensive computation of the unrestricted asymetric case.
\end{Ex}
\clearpage
\begin{Ex}
For \emph{Sidon-greedy} ($k=2$) we take $F = f = x_4+x_1-x_3-x_2$.  Thus $ac$ is simply $x_4+x_1= x_3+x_2$. Here, we get two distinct graphs for the involutions. The initial positions $(0,0),(1,1),(2,3),(3,2)$ are shared by all three variations. Then we get $(4,5),(5,9),(6,14)$ for the unrestricted, $(4,5),(5,4),(6,9)$ for max $ac$ and $(4,4),(5,7),(6,6)$ for order preserving respectively. Intuitively, the values of the order preserving greedy injection would tend to concentrate more than those of max $ac$. However, our computations so far show that max $ac$ paperers to be less ``dense".  Will such tendencies persist?
\end{Ex}
\clearpage
\begin{Ex}\label{ex4}
For \emph{line-greedy} we take $ac = e_1 \text{ OR } e_2 \text{ OR }\ldots $ where $$f_i(x_1,x_2,x_3) = \alpha_i (x_1-x_3)+\beta_i( x_2- x_3),$$ for an enumeration of all relatively prime $\alpha_i,\beta_i$. Another way to express this greedy algorithm is that $\pi_g(n)$ takes the least integer which does not lie on any line defined by a distinct pair of lattice points of the form $((i,\pi_g(i)), (j,\pi_g(j)))$, for $i,j\in \{0,1,\ldots ,n-1\}$. (Here max $ac$ and order preserving coincides.) Both sequences begin $$(0,0), (1,1), (2,3), (3,2), (4,5), (5,4), (6,8), (7,11), (8,6), (9,13), (10,12), (11,7).$$ The next entry differs, it is $(12,22)$ and $(12,10)$ respectively. Note that $(12,9)$ is impossible for both sequences since $(2,3)$ and $(6,8)$ are contained in both, but that $(12,10)$ is impossible for the unrestricted avoidance since $(9,13), (10,12),(12,10)$ is a decreasing progression. Although max $ac$ is a weaker than unrestricted $ac$, the respective values for $\pi_g$ seem to ``spread out" at roughly the same rate. Will this persist for greater values?   
\end{Ex} 
\clearpage
\begin{Ex}\label{ex5}
In \emph{parallel-greedy}, $\pi_g$ avoids any four lattice points that can be ordered pairwise as to define two parallel lines. This is a vast strengthening of the Sidon condition.
\end{Ex}
\clearpage

\begin{Thm}\label{thm42}
Suppose that an arithmetic condition $ac$ is given in $k$ variables and that $\pi_g$ is the greedy injection which avoids $\max ac$, or the order preserving $ac$. Then $\pi_g$ is an involution. 
\end{Thm}

\noindent{\bf Proof.} In this proof we write $\pi$ for $\pi_g$ with the max $ac$ condition. 
Suppose first that there is an $m<n$ such that 
\begin{align}\label{pimn}
\pi(m)=n,
\end{align}
and that $\pi(\pi(i))=i$ for each $i<n$ with $\pi(i) < \pi(m)$. 
We need to prove that $\pi(n)=m$. If this does not hold, by our assumption, we must have $\pi(n)>m$. This means that $\pi$ rejects $m$ for input $n$. Hence (\ref{x}), (\ref{pi}) and (\ref{<}) are satisfied simultaneously for some $j$, but where $\pi(n)$ is exchanged for $m$. Hence (\ref{pi}) becomes $$ac(\pi(x_{1}),\ldots,\pi(x_{{j}}),m,\pi(x_{{j+1}}),\ldots, \pi(x_{{k-1}})),$$ with $\pi(x_i)<m<n=\pi(m)$ for all such $x_i$, which by assumption implies $\pi(\pi(x_i)) = x_i$. Hence (\ref{x}) becomes $$ac(\pi(\pi(x_{1})),\ldots,\pi(\pi(x_{{j}})),n ,\pi(\pi(x_{{j+1}})),\ldots, \pi(\pi(x_{{k-1}})), $$ which contradicts the choice of $n$ for input $m$ in (\ref{pimn}).%

Suppose next, for a contradiction, that there is no $m$ such that (\ref{pimn}) holds, but there is an $m<n$ such that 
\begin{align}\label{nm}
\pi(n)= m. 
\end{align}
We assume that $\pi(\pi(i))=i$ for each $i<n$ with $\pi(i) < \pi(n)$. Then $\pi(m)>n$. Hence $n$ must have been rejected by $\pi$ for input $m$. That is, by (\ref{nm}), there must exist simultaneous solutions to  
$$ac(x_{1},\ldots,x_{{j}},m ,x_{{j+1}},\ldots, x_{{k-1}})=$$$$ac(\pi(\pi(x_{1})),\ldots,\pi(\pi(x_{{j}})),\pi(n) ,\pi(\pi(x_{{j+1}})),\ldots, \pi(\pi(x_{{k-1}}))) $$ 
and $$ac(\pi(x_{1}),\ldots,\pi(x_{{j}}),n,\pi(x_{{j+1}}),\ldots, \pi(x_{{k-1}})),$$ which contradicts the definition of $\max ac$. The proof for order preserving $ac$ is identical. \hfill $\Box$\\

\section{Discussion}
In \cite{La} it is demonstrated that symmetric extensions of Wythoff Nim have ``less dense" P-positions than those of Wythoff Nim. Here we have introduced an asymmetric variant of Wythoff Nim. The proof in \cite{La} follows an intuitive approach that we believe applies also for asymmetric cases (the final question in \cite{KnLa04}). Many combinatorial heap games (e.g. variations of Nim, Wythoff Nim) have symmetric rules and hence symmetric outcomes. It is interesting that Theorem \ref{thm42} emphasizes the relation of symmetry to that, in the game setting, heap sizes decrease for each move. 

One can also study our two heap games in a symmetric setting without any Nim type moves. We did not yet pursue any such problems.

Another theme of this paper has been to develop a framework for studying how the ``density" of P-positions (lattice points) changes when game rules are altered, by say heap-size invariant modifications. Combinatorial number theory may benefit from this approach by studying variations of known sequences that have natural settings among heap games. New combinatorial games have been developed to emulate classical sequences, suggesting new interactions between combinatorial games and combinatorial number theory. We note that the games in Section \ref{s6} can easily be modified to several dimensions and similarly for the greedy algorithm, the latter, which seems to be more of a technical challenge. Many problems for the CGT part of the work has been left for the future, such as an investigation of nim-values, partizan variants (producing four outcome classes), scoring play games and mis\`{e}re analogies for any setting.

\end{document}